%% file: article_cmcg_maxwell.tex
\documentclass{amsart}

\usepackage{a4wide}
\usepackage[foot]{amsaddr}

\numberwithin{figure}{subsection}

\usepackage{pgfplots}
\usepackage{subcaption}
\usepackage{graphicx}

\usepackage{algorithm}
\usepackage{algorithmic}

\usepackage{amssymb}
\usepackage{mathrsfs}

\usepackage{hyperref}

\title{A controllability method for Maxwell's equations}

\author{T. Chaumont-Frelet$^{\dagger,\ddagger}$}
\author{M.J. Grote$^{\sharp}$}
\author{S. Lanteri$^{\dagger,\ddagger}$}
\author{J.H. Tang$^{\flat}$}

\address{\vspace{-.5cm}}
\address{\noindent \tiny \textup{$^\dagger$Inria, 2004 Route des Lucioles, 06902 Valbonne, France}}
\address{\noindent \tiny \textup{$^\ddagger$Laboratoire J.A. Diedonn\'e, Parc Valrose, 28 Avenue Valrose, 06108 Nice Cedex 02, France}}
\address{\noindent \tiny \textup{$^\sharp$Department of Mathematics and Computer Science, University of Basel, Spiegelgasse 1, 4051 Basel, Switzerland}}
\address{\noindent \tiny \textup{$^\flat$ISTerre, 1381 Rue de la Piscine, 38610 Gières, France}}

\input{general_commands}

\input{specific_commands}

\begin{document}

\maketitle
\thispagestyle{empty}

\begin{abstract}
We propose a controllability method for the numerical solution of time-harmonic Maxwell's
equations in their first-order formulation. By minimizing a quadratic cost functional,
which measures the deviation from
periodicity, the controllability method determines iteratively a periodic
solution in the time domain. At each conjugate gradient
iteration, the gradient of the cost functional is simply computed by running any
time-dependent simulation code forward and backward for one period,
thus leading to a non-intrusive implementation easily integrated into existing software.
Moreover, the proposed algorithm automatically inherits
the parallelism, scalability, and low memory footprint of the underlying
time-domain solver. Since the time-periodic solution obtained by minimization
is not necessarily unique, we apply a cheap post-processing filtering
procedure which recovers the time-harmonic solution from any minimizer.
Finally, we present a series of numerical examples which show
that our algorithm greatly speeds up the convergence towards the desired time-harmonic
solution when compared to simply running the time-marching code until the time-harmonic
regime is eventually reached.

\vspace{.2cm}

\end{abstract}

\input{1-intro}
\input{2-description}
\input{3-model}
\input{4-abstract-analysis}
\input{5-controllability-method}
\input{6-numerical-examples}
\input{7-conclusion}

\bibliographystyle{amsplain}
\bibliography{bibliography.bib}
\end{document}

%% file: general_commands.tex
\renewcommand{\Re}{\operatorname{Re}}
\renewcommand{\Im}{\operatorname{Im}}

\newcommand{\eq}{:=}

\newcommand{\grad}{\boldsymbol \nabla}

\newcommand{\curl}{\grad \times}

\newcommand{\ccurl}{\boldsymbol{\operatorname{curl}}}

\newcommand{\jmp}[1]{\,[\![#1]\!]}
\newcommand{\avg}[1]{\{\!\!\{#1\}\!\!\}}

\newcommand{\BbBC}{\mathbb{C}}

\newcommand{\BbBR}{\mathbb{R}}

\newcommand{\BE}{\boldsymbol E}

\newcommand{\BG}{\boldsymbol G}
\newcommand{\BH}{\boldsymbol H}

\newcommand{\BJ}{\boldsymbol J}

\newcommand{\BL}{\boldsymbol L}

\newcommand{\BY}{\boldsymbol Y}
\newcommand{\BZ}{\boldsymbol Z}

\newcommand{\bc}{\boldsymbol c}
\newcommand{\bd}{\boldsymbol d}
\newcommand{\be}{\boldsymbol e}
\newcommand{\bg}{\boldsymbol g}
\newcommand{\bh}{\boldsymbol h}

\newcommand{\bj}{\boldsymbol j}

\newcommand{\bn}{\boldsymbol n}
\newcommand{\bo}{\boldsymbol o}
\newcommand{\bp}{\boldsymbol p}
\newcommand{\bq}{\boldsymbol q}

\newcommand{\bt}{\boldsymbol t}

\newcommand{\bv}{\boldsymbol v}
\newcommand{\bw}{\boldsymbol w}
\newcommand{\bx}{\boldsymbol x}

\newcommand{\CF}{\mathcal F}

\newcommand{\CH}{\mathcal H}

\newcommand{\CP}{\mathcal P}

\newcommand{\CT}{\mathcal T}

\newcommand{\CX}{\mathcal X}

\newcommand{\LA}{\mathscr A}

\newcommand{\LD}{\mathscr D}

\newcommand{\LF}{\mathscr F}
\newcommand{\LG}{\mathscr G}
\newcommand{\LH}{\mathscr H}

\newcommand{\LJ}{\mathscr J}
\newcommand{\LK}{\mathscr K}
\newcommand{\LL}{\mathscr L}

\newcommand{\LR}{\mathscr R}

\newcommand{\LU}{\mathscr U}
\newcommand{\LV}{\mathscr V}
\newcommand{\LW}{\mathscr W}

\newcommand{\LY}{\mathscr Y}
\newcommand{\LZ}{\mathscr Z}

\newcommand{\BCH}{\boldsymbol{\CH}}

\newcommand{\BCP}{\boldsymbol{\CP}}

\newcommand{\BCX}{\boldsymbol{\CX}}



%% file: specific_commands.tex
\newcommand{\ee}{\boldsymbol \varepsilon}
\newcommand{\mm}{\boldsymbol \mu}

\newcommand{\sig}{\boldsymbol \sigma}
\newcommand{\eps}{\varepsilon}

\newcommand{\xxi}{\boldsymbol \xi}

\newcommand{\GP}{{\Gamma_{\rm P}}}
\newcommand{\GI}{{\Gamma_{\rm I}}}

\newcommand{\CFPh}{{\CF_{{\rm P},h}}}
\newcommand{\CFIh}{{\CF_{{\rm I},h}}}

\newcommand{\CFih}{{\CF_{{\rm int},h}}}

\newcommand{\Cstab}{C_{\rm stab}}

\newtheorem{theorem}{Theorem}

\newtheorem{lemma}[theorem]{Lemma}
\newtheorem{assumption}[theorem]{Assumption}

\numberwithin{theorem}{section}
\numberwithin{equation}{section}

\newcommand{\scurl}{\operatorname{curl}}
\newcommand{\vcurl}{\boldsymbol{\scurl} \;}

\newcommand{\magn}{\bh}

\newcommand{\bphi}{\boldsymbol \phi}

\newcommand{\cgit}{\ell}    

\newcommand{\enorm}[1]{\left |\!\left |\!\left | #1 \right |\!\right |\!\right |}

\newcommand{\VUh}{{\mathbb U}_h}

\newcommand{\LVI}{\LV_{\rm I}}

\newcommand{\zero}{\boldsymbol 0}

%% file: 1-intro.tex
\section{Introduction}

Efficient numerical methods for
electromagnetic wave propagation are central to a wide range of
applications in science and technology \cite{assous_ciarlet_labrunie_2018a,griffiths_1999a}.
For wave phenomena with harmonic time dependence, governed by a single angular frequency $\omega > 0$,  the electromagnetic wave field satisfies time-harmonic Maxwell's equations in a
domain $\Omega \subset \mathbb R^3$: Given a current density $\bj: \Omega \to \mathbb C^3$,
we seek two vector fields $\be,\bh: \Omega \to \mathbb C^3$ such that
\begin{subequations}
\label{eq_maxwell_freq_strong}
\begin{equation}
\label{eq_maxwell_freq_strong_volume}
\left \{
\begin{array}{rcl}
i\omega \ee \be + \sig \be + \curl \bh &=& \bj,
\\
i\omega \mm \bh - \curl \be &=& \zero,
\end{array}
\right .
\end{equation}
inside the computational domain $\Omega$, where the first-order tensors $\ee$, $\sig$ and $\mm$
are the permittivity, conductivity and permeability
of the medium in $\Omega$. At the boundary $\partial \Omega$
of $\Omega$, divided into two disjoint sets $\GP$ and $\GI$, we impose
the boundary conditions
\begin{equation}
\label{eq_maxwell_freq_strong_boundary}
\left \{
\begin{array}{rcll}
\be \times \bn &=& \zero & \text{ on } \GP,
\\
\be \times \bn + \BZ \bh_\tau &=& \bg & \text{ on } \GI,
\end{array}
\right .
\end{equation}
\end{subequations}
where $\bn$ stands for the outward unit normal to $\partial \Omega$
and $\bh_\tau \eq \bn \times (\bh \times \bn)$. Here, the
first-order tensor $\BZ$, defined on $\GI$, describes a surface impedance
while $\bg: \GI \to \mathbb C^3$ typically represents
incident electromagnetic field. The PEC condition on $\GP$ corresponds to the surface of
a perfectly conducting material whereas the impedance
boundary condition on $\GI$ either models the boundary of an imperfect
conductor or corresponds to an approximation of the Silver-M\"uller radiation
condition \cite{colton_kress_2012a}. Note that $\GP$ or $\GI$ may be empty.
%

In heterogeneous media with intricate
geometries, Galerkin discretizations based on variational formulations of
\eqref{eq_maxwell_freq_strong}, such as curl-conforming finite elements or
discontinous Galerkin (DG) methods \cite{li_lanteri_perrussel_2014a,monk_2003a},
probably are the most flexible and competitive approaches
currently available. If $\omega$ is ``large'' and
the computational domain spans many wavelengths,
resolving the wavelength and limiting dispersion errors requires the
use of highly refined meshes coupled with high-order elements
\cite{chaumontfrelet_nicaise_2019a,melenk_sauter_2020a}. Hence, the high-frequency regime
typically leads to large, sparse, indefinite and ill-conditioned linear systems which need to be solved numerically
by direct or iterative methods.
Although considerable progress has been achieved over the past decades
\cite{amestoy_ashcraft_boiteau_buttari_lexcellent_weisbecker_2015a,amestoy_duff_lexcellent_2000a},
the parallel implementation of scalable direct solvers remains a challenge when the number of
unknowns is large. On the other hand, the design of robust and efficient preconditioners for
iterative solvers is a delicate task \cite{ernst_gander_2012a}. Recent developments
include domain decomposition \cite{bonazzoli_dolean_graham_spence_tournier_2019a,li_lanteri_perrussel_2014a},
shifted-laplacian \cite{gander_graham_spence_2015a}, and sweeping \cite{tsuji_engquist_ying_2012a}
preconditioners. Still, the efficient solution of 3D time-harmonic Maxwell's equations with hetereogeneous
coefficients remains to this day a formidable challenge, especially in the high-frequency regime.

To avoid these difficulties, we instead transform \eqref{eq_maxwell_freq_strong} back to the time-domain and
consider its time-dependent counterpart
\begin{equation}
\label{eq_maxwell_time_strong}
\left \{
\begin{array}{rcll}
\ee \dot \BE + \sig \BE + \curl \BH &=& \BJ & \text{ in } \mathbb R_+ \times \Omega,
\\
\mm \dot \BH - \curl \BE &=& \zero & \text{ in } \mathbb R_+ \times \Omega,
\\
\BE \times \bn &=& \zero & \text{ on } \mathbb R_+ \times \GP,
\\
\BE \times \bn + \BZ \BH_\tau &=& \BG &\text{ on } \mathbb R_+ \times \GI,
\end{array}
\right .
\end{equation}
with time-harmonic forcing
$\BJ(t,\bx) \eq \Re \left \{\bj(\bx) e^{i\omega t}\right \}$,
$\BG(t,\bx) \eq \Re \left \{\bg(\bx) e^{i\omega t}\right \}$,
and initial conditions $\BE|_{t=0} = \BE_0$ and $\BH|_{t=0} = \BH_0$ yet to
be specified. The key advantage of this strategy is that it only requires the
solution of a time evolution problem for which efficient numerical schemes, such as
finite differences \cite{taflove_hagness_2005a,yee_1966a} or DG
\cite{fezoui_lanteri_lohrengel_piperno_2005a,grote_schneebeli_schotzau_2007a,hesthaven_warburton_2002a}
discretizations coupled with explicit time integration, can be utilized.
As these algorithms are inherently parallel with a low memory footprint,
they are extremely attractive on modern computer architectures.

In this context, a simple and common approach follows from the limiting amplitude principle
\cite{morawetz_1962a}, which states under suitable assumptions that the solution of \eqref{eq_maxwell_time_strong} ``converges'' to the time-harmonic solution
in the sense that $\BE(t,\bx) \to \Re \left \{\be(\bx)e^{i\omega t}\right \}$ and
$\BH(t,\bx) \to \Re \left \{\bh(\bx)e^{i\omega t}\right \}$
as $t \to +\infty$. Thus, to solve  \eqref{eq_maxwell_freq_strong} one can simply
simulate time-dependent Maxwell's equations
for a ``sufficiently long'' time and eventually extract the time-harmonic solution.
However, as the final simulation time
required to obtain
an accurate approximation may be very large, especially near resonances or in the presence of trapping geometries, the usefulness of this approach is somewhat limited \cite{bardos_rauch_1994a}.

Both controllability methods and fixed-point iterations
have been proposed to accelerate convergence and determine initial conditions ($\BE_0,\BH_0$)
which render the time-dependent solution $T$-periodic with period $T \eq 2\pi/\omega$.
Inspired by the seminal work in \cite{lions_1988a}, controllability
methods (CM) \cite{bristeau_glowinski_periaux_1994a,bristeau_glowinski_periaux_1998a} reformulate
the controllability problem as a minimization problem
for a quadratic cost functional $J(\BE_0,\BH_0)$, which measures the
misfit between $(\BE_0,\BH_0)$ and the time-dependent solution $(\BE(T),\BH(T))$
after one period. Then, the functional $J$ is minimized
by a conjugate gradient (CG) iteration, which leads to the combined
controllability method-CG algorithm, or CMCG for short. Alternatively, fixed-point iterations
determine the $T$-periodic solution by applying a judicious filtering operator at
each iteration to achieve convergence  \cite{AppeloGarciaRunborg, peng2021emwaveholtz}. As the convergence of fixed-point iterations can be slow
near resonances or in the presence of trapping geometries, an outer CG or GMRES Krylov subspace method
must be applied, depending on boundary condititions.

When using the controllability approach, one faces two central questions:
efficient computation of the gradient $J'$ and uniqueness of the time-periodic solution.
As early work on CMCG methods was restricted to scattering problems
from acoustics \cite{bristeau_glowinski_periaux_1994a,bristeau_glowinski_periaux_1998a}
or electromagnetics \cite{bristeau:inria-00073072} in second-order formulation, the computation
of  $J'$ always required the solution of a strongly elliptic (coercive) problem.
To avoid solving that additional
elliptic problem at each CG iteration, the controllability method was later applied to the Helmholtz
equation in first-order formulation \cite{kahkonen_glowinski_rossi_makinen_2011a} using
Raviart-Thomas FE for the spatial discretization; due to the lack of available mass-lumping,
however, the mass-matrix then needed to be inverted at each time-step during the
time integration. By combining a first-order formulation with a DG discretization, a scalable
parallel formulation was recently derived \cite{grote_nataf_tang_tournier_2020a},
which completely avoids the need for solving any elliptic problem or inverting the mass-matrix.

In general, the $T$-periodic solution of \eqref{eq_maxwell_time_strong} is not unique
and hence does not necessarily yield the desired (unique) time-harmonic solution  of
\eqref{eq_maxwell_freq_strong}. For sound-soft acoustic scattering, where Dirichlet
and impedance conditions are imposed on distinct parts of the boundary, the $T$-periodic
solution in fact is unique and the one-to-one correspondence is therefore immediate.
For other boundary-value problems, however,
such as sound-hard scattering or problems in bounded physical domains,
the periodic solution is generally no longer unique, as it may contain additional
($T$-periodic) spurious modes. Two ideas have been
proposed as a remedy to extend the CMCG approach to arbitrary boundary conditions.
First, uniqueness can be restored by modifying $J$, though at a small price in
the computation of its gradient \cite{bardos_rauch_1994a,grote_tang_2019a}. Alternatively,
a cheap filtering operator can be applied as a post-processing step to any minimizer of $J$,
which removes any spurious modes \cite{grote_nataf_tang_tournier_2020a,tang_2020a}
and thus restores uniqueness using the original cost functional $J$.

Here we propose a CMCG method for time-harmonic
Maxwell's equations  \eqref{eq_maxwell_freq_strong} in their first order formulation,
which completely avoids the solution of any elliptic problem, and combine it with a
post-processing filtering step to guarantee uniqueness, regardless of the boundary conditions.
Moreover, thanks to a DG discretization in space, the mass-matrix is automatically block-diagonal.
Hence, the resulting CMCG algorithm is inherently parallel and scalable but also
guaranteed to converge to the time-harmonic solution
starting from any initial guess, as long as time-harmonic Maxwell's equations
\eqref{eq_maxwell_freq_strong} are well-posed for the frequency $\omega$ under
consideration.

The remainder of this work is organized as follows. We provide a formal
description of the algorithm and a discussion of our key theoretical
results in Section \ref{section_description}. As the mathematical framework required to
rigorously define and analyze Maxwell's equations is rather involved,
the precise description and preliminary results are postponed to Section \ref{section_model}.
Section \ref{section_perio} contains the bulk of the theory, where we carefully
analyze the relation between the time-harmonic and time-periodic solutions.
Here, our contributions are twofold. On the one hand, we identify configurations
of boundary conditions and right-hand sides for which the unique time-periodic solution
coincides with the time-harmonic solution. On the other hand, we show that
the filtering procedure introduced in \cite{grote_nataf_tang_tournier_2020a,tang_2020a}
always recovers the time-harmonic solution from any minimizer, as long as
\eqref{eq_maxwell_freq_strong} is well-posed.
In Section \ref{section_CMCG}, we describe in detail our CMCG method and establish
its convergence toward the time-harmonic solution.
In Section \ref{section_numer},
we present various numerical experiments highlighting the performance of the proposed CMCG algorithm.
Here, we  benchmark the proposed CMCG algorithm against the limiting amplitude principle,
where pure time-marching (without controllability) is utilized, as both methods are
non-invasive and easily integrated with any existing time-marching code; in contrast, efficient preconditioners typically require an important and dedicated implementation effort. Finally, we provide
in Section \ref{section_conclusion} some concluding remarks.

%% file: 2-description.tex
\section{Main results}
\label{section_description}

Throughout this work, we adopt the notation $U = (\be,\bh)$ for a time-harmonic
electromagnetic field, while the calligraphic font $\LU = (\BE,\BH)$ is reserved
for time-dependent fields. It is easily seen that if $U$ is a time-harmonic
field solution to \eqref{eq_maxwell_freq_strong} with right-hand side $\bj$ and $\bg$, then
$\LU(t,\bx) \eq \Re \{U(\bx) e^{i\omega t}\}$ is the solution of time-dependent
Maxwell's equations \eqref{eq_maxwell_time_strong} with right-hand side
$\BJ(t,\bx) \eq \Re \{\bj(\bx) e^{i\omega t}\}$, $\BG(t,\bx) \eq \Re \{\bg(\bx) e^{i\omega t}\}$,
and initial condition $\LU_0 \eq \Re U$.

The CMCG algorithm hinges on an idea that is essentially the converse of the above
statement. Namely, we seek an initial condition $\LU_0$ such that the resulting time-dependent
field $\LU$ (with right-hand sides $\BJ$ and $\BG$ as above) is time-periodic, with period
$T \eq 2\pi/\omega$.
Let $P_{\bj,\bg,\omega}: \LU_0 \to \LU(T)$ denote the (affine) operator mapping the initial
condition $\LU_0$ to the solution $\LU$ of \eqref{eq_maxwell_time_strong} with
time-harmonic right-hand sides $\BJ$ and $\BG$ evaluated at time $T$. Then, the
``controllability method'' corresponds to solving (linear) equation
$P_{\bj,\bg,\omega} \LU_0 = \LU_0$.

At this point, three main questions arise. First, if the time-dependent solution with
initial condition $\LU_0$ is periodic, can we ensure that $\LU_0 = \Re U$, where $U$
is the corresponding frequency-domain solution?  Second, can we design an efficient
algorithm to solve for $P_{\bj,\bg,\omega} \LU_0 = \LU_0$? Finally, can we prove the
convergence of this algorithm?

\subsection{The structure of periodic solutions}

Our first set of results characterizes those initial conditions $\LU_0$ such that
$\LU_0 = P_{\bj,\bg,\omega} \LU_0$. In essence, we establish that
\begin{equation*}
\LU_0 = \Re \left ([\bp,\bq] + U + \sum_{|\ell| \geq 2} U_\ell\right ),
\end{equation*}
where $U$ is the unique time-harmonic solution, $\bp$ and $\bq$ are two curl-free fields with
$\bp \times \bn = \bq \times \bn =0$ on $\GI$, and for all $|\ell| \geq 2$, $U_\ell$ is any
time-harmonic solution with frequency $\ell\omega$ and vanishing right-hand sides.
Thus, if time-harmonic problem \eqref{eq_maxwell_freq_strong} is well-posed for
all multiples $\ell\omega$ of $\omega$, then we simply have
$\LU_0 = \Re \left ([\bp,\bq] + U\right )$,
which holds whenever the problem features dissipation ($\operatorname{supp} \sig \neq \emptyset$
and/or $|\GI| > 0$). Moreover, we show that if both $\LU_0$ and $\bj$ are
orthogonal to curl-free fields, then $\bp = \bq = \bo$, so that $\LU_0 = \Re U$.
In fact, if $\Omega$ is simply connected, we have $\bp = \grad p$
and $\bq = \grad q$ for two scalar functions $p$ and $q$, while the condition on
$\LU_0$ and $\bj$ simply means that they are divergence-free.

Our second set of results concerns the post-processing of periodic solutions
by the filtering operator
\begin{equation}
\label{eq_filter_rhs}
F_{\bj,\bg,\omega} \LU_0 \eq \frac{2}{T} \int_0^T \LU(t) e^{-i \omega t} dt,
\end{equation}
where $\LU$ is the solution to time-dependent Maxwell's equations \eqref{eq_maxwell_time_strong}
with initial condition $\LU_0$ and right-hand sides $\BJ$ and $\BG$. Note that $F_{\bj,\bg,\omega}$
may be easily computed ``on the fly'' during time-marching while computing $P_{\bj,\bg,\omega}$
without storing the time-history of $\LU(t)$. Then, our key result states that
$U = F_{\bj,\bg,\omega} \LU_0$ for any initial condition $\LU_0$ satisfying
$\LU_0 = P_{\bj,\bg,\omega} \LU_0$, as long as time-harmonic problem \eqref{eq_maxwell_freq_strong}
is well-posed for the frequency $\omega$,

In fact, we prove the slightly stronger result that for any
initial condition $\LU_0$, $F_{\bj,\bg,\omega} \LU_0$ solves time-harmonic Maxwell's
equations with a modified right-hand side, where the misfit $(I-P_{\bj,\bg,\omega})\LU_0$ is added
to the physical source terms. This result enables us to control the error
$U - F_{\bj,\bg,\omega} \LU_0$ by the misfit $\LU_0 - P_{\bj,\bg,\omega} \LU_0$.
It is also central for subsequently analyzing the convexity of the cost functional.

\subsection{The CMCG algorithm}
\label{sec_CMCG_algo}

To determine an initial condition $\LU_0$ that leads to a time-periodic solution,
i.e. $\LU_0=P_{\bj,\bg,\omega}\LU_0$, we minimize the ``energy functional''
\begin{equation*}
J(\LU_0)
\eq
\frac{1}{2} \|\LU(T)-\LU_0\|_{\ee,\mm}^2
=
\frac{1}{2} \|(I-P_{\bj,\bg,\omega})\LU_0\|^2_{\ee,\mm}
\end{equation*}
which measures the ($\ee,\mm$-weighted) $L^2(\Omega)$-misfit between the initial condition
and the solution after one period. Since $P_{\bj,\bg,\omega}$ is an affine operator,
it can be decomposed as $P_{\bj,\bg,\omega} \LU_0 = P_\omega \LU_0 + \LG$, where
$\LG \eq P_{\bj,\bg,\omega} 0$ and the operator
$P_\omega \eq P_{\boldsymbol{0},\boldsymbol{0},\omega}$, which corresponds to the propagation
of the initial condition $\LU_0$ a time $T$ with zero right-hand side, is now linear.
Hence
\begin{equation*}
J(\LU_0)
=
\frac{1}{2} \|(I-P_\omega) \LU_0 - \LG\|_{\ee,\mm}^2,
\end{equation*}
is a standard quadratic functional.

The gradient is given by
\begin{equation*}
J'(\LU_0) = (I-P_\omega^\star)(I-P_\omega) \LU_0 - \LG^\star,
\qquad
\LG^\star \eq (I-P_\omega^\star) \LG,
\end{equation*}
where $P_\omega^\star$ denotes the adjoint of $P_\omega$, which
actually maps the final condition $\LW_T$ to $\LW(0)$ by back-propagation.
In practice the action of $P_\omega$ and $P_\omega^\star$ on any $\LU_0$ is
simply obtained by solving \eqref{eq_maxwell_time_strong} numerically in the
time-domain for one period. Hence, after the initialization
step described in Algorithm \ref{al_init}, we simply compute the gradient of $J$ by one
forward and one backward solve as listed in Algorithm \ref{al_grad}.

Once we have an efficient algorithm to compute $J'$, we may choose any quadratic minimization
algorithm \cite{ciarlet_1989a}. Here, we employ the conjugate gradient method,
resulting in Algorithm \ref{al_cmcg}. Note that in practice the evaluation
of the scalar product $(\LU_0,\LV_0)_{\ee,\mm}$ simply amounts to computing
$\mathbb V^{\rm T} \mathbb M \mathbb U$, where $\mathbb M$ is the mass matrix
arising from space discretization, and $\mathbb U$ (resp. $\mathbb V$) is
the discrete vector of degrees of freedom representing $\LU_0$ (resp. $\LV_0$).

\begin{algorithm}[t]
\begin{algorithmic}[1]
\REQUIRE right-hand sides $\bj$ and $\bg$
\STATE compute $\LG = P_{\bj,\bg,\omega} 0$ by time-marching for one period
\STATE compute $\LG_T = P_\omega^\star \LG$ by back-propagating over one peroid
\STATE set $\LG^\star = \LG-\LG_T$
\RETURN $\LG^\star$
\end{algorithmic}
\caption{Initialization}
\label{al_init}
\end{algorithm}
\begin{algorithm}[t]
\begin{algorithmic}[1]
\REQUIRE real-valued eletromagnetic field $\LU_0$, precomputed $\LG^\star$
\STATE compute $\LU_T = P_\omega \LU$ by time-marching for one period
\STATE set $\LW_T = \LU_T-\LU_0$.
\STATE compute $\LW_0 = P_\omega^\star \LW_T$ by back-propagation over one period
\STATE set $J'(\LU_0) = \LW_T-\LU_0-\LG^\star$.
\RETURN $J'(\LU_0)$
\end{algorithmic}
\caption{Gradient evaluation}
\label{al_grad}
\end{algorithm}
\begin{algorithm}[t]
\caption{CMCG Algorithm}
\label{al_cmcg}
\begin{algorithmic}[1]
\REQUIRE
right-hand sides $\bj$ and $\bg$,
initial guess $\LU_0^{(0)}$,
tolerance $\delta$,
maximum iteration $\cgit_{\rm max}$
\STATE compute $\LG^\star$ from $\bj$ and $\bg$ with Algorithm \ref{al_init}
\STATE compute $\LJ' = J'(\LU_0^{(0)})$ with Algorithm \ref{al_grad}
\STATE set $\LR^{(0)} = \LJ'$, $\LD^{(0)} = \LJ'$
\FOR{$\cgit=0,\dots,\cgit_{\rm max}-1$}
\IF{$\|\LR^{(\cgit)}\|_{\ee,\mm} \leq \delta \|\LR^{(0)}\|_{\ee,\mm}$}
\RETURN $\LU_0^{(\cgit)}$
\ENDIF
\STATE compute $\LA = J'(\LD^{(\cgit)})+\LG^\star$ with Algorithm \ref{al_grad}
\STATE set $\alpha = \|\LR^{(\cgit)}\|_{\ee,\mm}^2 / (\LD^{(\cgit)},\LA)_{\ee,\mm}$
\STATE set $\LU_0^{(\cgit+1)} = \LU_0^{(\cgit)} + \alpha \LD^{(\cgit)}$
\STATE set $\LR^{(\cgit+1)} = \LR^{(\cgit)} - \alpha \LA$
\STATE set $\beta = \|\LR^{(\cgit+1)}\|_{\ee,\mm}^2/\|\LR^{(\cgit)}\|_{\ee,\mm}^2$
\STATE set $\LD^{(\cgit+1)} = \LR^{(\cgit)} + \beta \LD^{(\cgit)}$
\ENDFOR
\RETURN $\LU_0^{(\cgit_{\rm max})}$
\end{algorithmic}
\end{algorithm}

\subsection{Convexity of the functional and convergence}

Finally, we address the convexity of the energy functional, which immediately
relates to the convergence of the CMCG algorithm. It has been previously established
that $J$ is strongly convex for the case of sound-soft scattering by a convex obstacle,
but that it is \emph{not} necessarily so for general geometries \cite{bardos_rauch_1994a}.
Here, we show that $J$ is strongly convex in an appropriate sense
as long as time-harmonic problem \eqref{eq_maxwell_freq_strong} is well-posed, thereby
ensuring the convergence of the proposed algorithm. To do so, we introduce
a second filtering operator $F_\omega \LU_0 \eq F_{\zero,\zero,\omega} \LU_0$
that is defined as \eqref{eq_filter_rhs}, but with right-hand sides $\bj = \bg = \zero$.
Our key result is that $J$ is continuous, uniformly-Lipschitz
and strictly convex on the space of initial conditions modulo the kernel of $F_\omega$.
This quotient space is only used as a technical tool in the proofs, and, in practice,
if $\LU_0^{(\cgit)}$ is the initial condition at iteration $\cgit$ in the CG algorithm, then
$F_{\bj,\bg,\omega} \LU_0^{(\cgit)} \to U$ for any initial guess $\LU_0^{(0)}$.

%% file: 3-model.tex
\section{Settings and preliminary results}
\label{section_model}

This section provides the mathematical framework needed to
rigorously anayze the CMCG algorithm.

\subsection{Domain and coefficients}

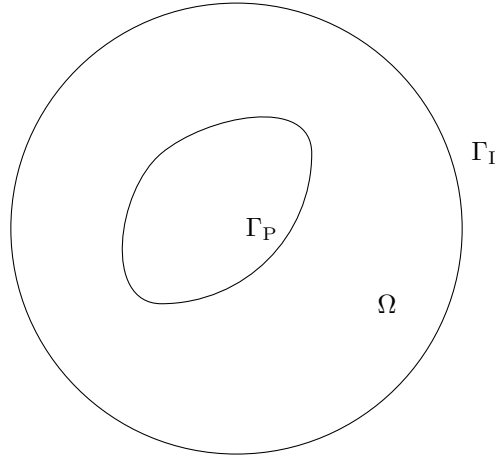
\begin{figure}
\begin{center}
\input{figures/boundary_condition_settings}
\end{center}
\caption{Example of boundary condition settings}
\label{figure_boundary}
\end{figure}

We consider time-harmonic Maxwell's equations set in a Lipschitz domain
$\Omega \subset \BbBR^3$. The boundary $\Gamma \eq \partial \Omega$ of $\Omega$
is partitioned into two relatively open disjoint subsets $\GP$ and $\GI$.
We assume that $\overline{\GP} \cap \overline{\GI} = \emptyset$,
which is not mandatory, but simplifies the analysis. Figure \ref{figure_boundary} presents
a possible configuration.

To avoid the proliferation of necessary notation to handle both two
and three-dimensional problems at the same time, we restrict our theoretical
investigations to three-dimensional domains However, our analysis also applies
to two-dimensional problems in any polarization with natural modifications.
For the sake of simplicity, we also avoid dealing with boundary sources
in our theoretical analysis, and focus on volumic sources.
Still, our numerical experiments show, that our CMCG method applies equally well
with both types of sources.

We consider three measurable symmetric tensor-valued functions
$\ee,\mm,\sig: \Omega \to \mathbb S(\BbBR^3)$ which respectively represent
the electric permittivity, the magnetic permeability, and the conductivity of
the material contained in $\Omega$. These tensors are assumed to be uniformly
bounded. We require that $\ee$ and $\mm$ are uniformly elliptic in $\Omega$.
For the conductivity, we assume that $\sig = \zero$ outside some
set $\Omega_{\sig} \subset \Omega$ with Lipschitz boundary
$\Gamma_{\sig} \eq \partial \Omega_{\sig}$ with $\sig$ uniformly elliptic
in $\Omega_{\sig}$.

On $\GI$, we consider a symmetric tensor-valued ``impedance'' function
$\BZ: \GI \to \mathbb S(\BbBR^3)$ which is assumed to be measurable with
respect to the surface measure, uniformly bounded and elliptic. We
also assume that $\BZ$ is tangential, i.e., for all $\xxi \in \BbBR^3$
and a.e. $\bx \in \GI$, $\xxi \cdot \bn(\bx) = 0$ implies that $\BZ(\bx) \cdot \xxi = 0$.
Finally, $\BY \eq \BZ^{-1}$ denotes the inverse of $\BZ$.

\subsection{Functional spaces}
\label{section_model_functional}

If $\mathbb K = \BbBR$ or $\BbBC$, $L^2(\Omega,\mathbb K)$
denotes the space of measurable square integrable functions
mapping $\Omega$ to $\mathbb K$ \cite{adams_fournier_2003a}.
Similarly, $L^2(\GI,\mathbb K)$ is the space of functions from
$\GI$ to $\mathbb K$ that are square integrable with respect to
the surface measure of $\GI$. For vector-valued function, we write
$\BL^2(\Omega,\mathbb K) \eq \left (L^2(\Omega,\mathbb K)\right )^3$
and
$\BL^2(\GI,\mathbb K) \eq \left (L^2(\GI,\mathbb K)\right )^3$.
We denote by $(\cdot,\cdot)_\Omega$ and $(\cdot,\cdot)_\GI$
the inner-products of these spaces. If $\bphi$ is a measurable essentially
bounded tensor, we employ the notations
$\|\cdot\|_{\bphi,\Omega}^2 = (\bphi\cdot,\cdot)_\Omega$
and $\|\cdot\|_{\bphi,\GI}^2 = (\bphi \cdot,\cdot)_\GI$.
As usual, $H^1(\Omega)$ stands for the first-order Sobolev space
\cite{adams_fournier_2003a}. If $\gamma \subset \partial \Omega$ is a relatively
open subset, $H^1_\gamma(\Omega,\mathbb K)$ is the subset of functions of
$H^1(\Omega,\mathbb K)$ with vanishing trace on $\gamma$.

For the analysis, we also need Sobolev spaces of vector-valued functions with ``well-defined''
curl, denoted by
$\BCH(\ccurl,\Omega,\mathbb K)
\eq
\left \{
\bv \in \BL^2(\Omega,\mathbb K)
\; | \;
\curl \bv \in \BL^2(\Omega,\mathbb K)
\right \}$,
see \cite{girault_raviart_1986a}. Following \cite{fernandes_gilardi_1997a}, we can define
the tangential trace of a function $\bv \in \BCH(\ccurl,\Omega,\mathbb K)$
on $\GP$ and $\GI$, and introduce
$\BCX(\Omega,\mathbb K)
\eq
\left \{
\bv \in \BCH(\ccurl,\Omega,\mathbb K)
\; | \;
\bv_\tau|_\GI \in \BL^2(\GI,\mathbb K)
\right \}$
and
$\BCX_\GP(\Omega,\mathbb K)
\eq
\left \{
\bv \in \BCX(\Omega,\mathbb K)
\; | \;
\bv_\tau|_\GP = \zero
\right \}$.

To simplify the discussion below, we finally introduce the product spaces
$L(\Omega) \eq \BL^2(\Omega,\BbBC) \times \BL^2(\Omega,\BbBC)$,
$\LL(\Omega) \eq \BL^2(\Omega,\BbBR) \times \BL^2(\Omega,\BbBR)$,
$V(\Omega) \eq \BCX_\GP(\Omega,\BbBC) \times \BCX(\Omega,\BbBC)$
and $\LV(\Omega) \eq \BCX_\GP(\Omega,\BbBR) \times \BCX(\Omega,\BbBR)$.
In the remaining of this work, we follow the convention introduced
above: if $Y(\Omega)$ is a space of complex-valued
electromagnetic fields, $\LY(\Omega)$ always denotes its real-valued counterpart.

The spaces $L$ and $\LL$ are equipped with the inner product
\begin{equation}
\label{eq_innerproduct}
([\bv,\bw],[\bv',\bw'])_{\ee,\mm}
\eq
(\ee\bv,\bv')_\Omega + (\mm\bw,\bw')_\Omega
\end{equation}
for all $[\bv,\bw],[\bv',\bw'] \in L(\Omega)$
and the associated norm $\|\cdot\|_{\ee,\mm}^2 = (\cdot,\cdot)_{\ee,\mm}$,
while we introduce the energy norm
\begin{align}
\label{eq_energy_norm}
\enorm{[\bv,\bw]}^2
&\eq
\omega^2\|\bv\|_{\ee,\Omega}^2
+
\|\bv_\tau\|_{\BY,\GI}^2
+
\|\curl \bv\|_{\mm^{-1},\Omega}^2
+
\|\sig\bv\|_{\ee^{-1},\Omega}^2
\\
\nonumber
&\;+
\omega^2 \|\bw\|_{\mm,\Omega}^2
+
\|\bw_\tau\|_{\BZ,\GI}^2
+
\|\curl \bh\|_{\ee^{-1},\Omega}^2
\end{align}
for all $[\bv,\bw] \in V(\Omega)$. We also introduce the subspace
\begin{equation*}
\LVI(\Omega)
\eq
\left \{
[\be,\bh] \in \LV(\Omega)
\; | \;
\be \times \bn + \BZ \bh_\tau = \zero
\;
\text{ on } \GI
\right \},
\end{equation*}
of fields satisfying impedance condition \eqref{eq_maxwell_freq_strong_boundary} on $\GI$.

Finally, if $\LY(\Omega)$ is any of the aforementioned real-valued
spaces, then $C^0(0,T;\LY(\Omega))$ and $C^1(0,T;\LY(\Omega))$ contain
functions from $[0,T]$ to $\LY(\Omega)$.

\subsection{Variational formulation}

We introduce the sesquilinear form
\begin{equation}
\label{eq_weak_maxwell}
a([\be,\bh],[\bv,\bw])
\eq
(\sig \be,\bv)
+
(\BY \be_\tau,\bv_\tau)_\GI
+
(\BZ \bh_\tau,\bw_\tau)_\GI
+
(\bh,\curl \bv)
-
(\be,\curl \bw)
\end{equation}
for all $[\be,\magn],[\bv,\bw] \in V(\Omega)$. Then,
the weak formulation of \eqref{eq_maxwell_freq_strong} is:
Find $[\be,\magn] \in V(\Omega)$ such that
\begin{equation*}
i\omega ([\be,\bh],[\bv,\bw]) + a([\be,\bh],[\bv,\bw])
=
(\bj,\bv) + (\BY \bg \times \bn,\bv_\tau)_\GI + (\BZ \bg,\bw_\tau)_\GI
\end{equation*}
for all $[\bv,\bw] \in V(\Omega)$. By using integration by parts,
we easily verify that
\begin{equation}
\label{eq_adjoint}
a([\bv,\bw],[\be,\bh]) = \overline{a([\be,-\bh],[\bv,-\bw])}
\end{equation}
for all $[\bv,\bw],[\be,\bh] \in V_{\rm I}(\Omega)$.

\subsection{Well-posedness}

Throughout this work, we assume that the time-harmonic
problem under consideration is well-posed for the chosen angular
frequency $\omega$.

\begin{assumption}[Well-posedness]
\label{eq_assumption_well_posedness}
For all $\phi \in L(\Omega)$, there exists a unique
$S_\omega \phi \in V(\Omega)$ such that
\begin{equation}
\label{eq_maxwell_weak}
i\omega (S_\omega \phi,w)_{\ee,\mm} + a(S_\omega \phi,w)
=
(\phi,w)_{\ee,\mm}
\quad \forall w \in V(\Omega).
\end{equation}
In addition, the stability estimate
\begin{equation}
\label{eq_frequency_stability}
\enorm{S_\omega \phi} \leq \Cstab \|\phi\|_{\ee,\mm}
\end{equation}
holds true.
\end{assumption}

In \eqref{eq_frequency_stability}, $\Cstab$ is a dimensionless
constant that depends on the frequency $\omega$, the shape
of the boundaries $\GP$ and $\GI$, and the physical coefficients $\ee$,
$\mm$ and $\sig$. Unless the entire domain contains a conductive materials
(i.e. $\Omega_{\sig} = \Omega$), the stability constant will increase with
the frequency. In the most favorable case of a non-trapping configuration
\cite{hiptmair_moiola_perugia_2010a,moiola_spence_2019a}, we have
\begin{equation*}
\Cstab \simeq \frac{\omega d_\Omega}{c},
\end{equation*}
where $c \eq 1/\sqrt{\varepsilon_{\max}\mu_{\max}}$ is the (minimal) wavespeed
and $d_\Omega$ is the diameter of the computational domain. If
$\lambda \eq c/\omega$ denotes the wavelength, $\Cstab$ is actually proportional
to the number of wavelengths $N_\lambda \eq d_\Omega/\lambda$
across $\Omega$. The stability constant can however
exhibit ``arbitrarily bad'' behaviour in more complicated geometries
(close to a resonance frequency when $\Omega_{\sig} \eq \emptyset$ and $\GI \eq \empty \emptyset$
for instance). We also mention that when considering two-dimensional geometries,
the two possible polarizations are equivalent to scalar Helmholtz problems, for
which a vast body of literature is now available (see, e.g., \cite{graham_pembery_spence_2019a}
and the references therein).

For future references, we note that the ``converse'' estimate to
\eqref{eq_frequency_stability}, namely
\begin{equation}
\label{eq_frequency_continuity}
\|\phi\|_{\ee,\mm} \leq \enorm{S_\omega \phi},
\end{equation}
holds true, as can be seen from the strong form of
time-harmonic Maxwell's equations \eqref{eq_maxwell_freq_strong} and
definition \eqref{eq_energy_norm} of the energy norm.

We finally observe that in view of \eqref{eq_adjoint}, the operator $S_\omega^\star$
defined for all $\phi \in L(\Omega)$ by the variational equation
\begin{equation*}
i\omega(w,S_\omega^\star \phi)_{\ee,\mm} + a(w,S_\omega^\star \phi)
=
(w,\phi)_{\ee,\mm}
\quad \forall w \in L(\Omega),
\end{equation*}
has a very similar structure to $S_\Omega$. In particular, \eqref{eq_frequency_stability}
and \eqref{eq_frequency_continuity} hold true for $S_\omega^\star$ too.

\subsection{Time-harmonic solution}

Henceforth, we consider a fixed right-hand side
$\psi \in L(\Omega)$, and denote by $U \in V(\Omega)$ the
associated solution satisfying
\begin{equation}
\label{eq_frequency_solution}
i\omega(U,w) + a(U,w) = (\psi,w)_{\ee,\mm} \quad \forall w \in V(\Omega),
\end{equation}
whose existence and uniqueness follows from Assumption \ref{eq_assumption_well_posedness}.

\subsection{Time-dependent solutions}
\label{section_model_time}

Although existence and uniqueness results for the time-dependent Maxwell's equations
\eqref{eq_maxwell_time_strong} are fairly standard, we provide some detail here,
since the final controllability method seeks an initial condition lying only in
the space $\LL(\Omega)$, so that solutions to \eqref{eq_maxwell_time_strong} can only
be defined in a very weak sense.

Following Sections 4.3.1 and 5.2.4 of \cite{assous_ciarlet_labrunie_2018a},
we introduce the unbounded operator
\begin{equation*}
A:
\LVI(\Omega) \ni [\be,\bh]
\to
[\ee^{-1}\sig\be+\ee^{-1} \curl \bh,-\mm^{-1} \curl \be] \in \LL(\Omega).
\end{equation*}
Then Hille-Yosida's theorem \cite[Theorem 4.3.2]{assous_ciarlet_labrunie_2018a} shows
that for all $\LU_0 \in \LVI(\Omega)$ and $\LF \in C^1(0,T,\LL(\Omega))$,
there exists a unique $\LU \in C^1(0,T,\LL(\Omega)) \cap C^0(0,T,\LVI(\Omega))$
such that
\begin{equation}
\label{eq_ode_maxwell_time}
\left \{
\begin{array}{rcll}
\dot \LU(t) + A\LU(t) &=& \LF(t) & t \in [0,T],
\\
\LU(0) &=& \LU_0, &
\end{array}
\right .
\end{equation}
and the estimate
\begin{equation}
\label{eq_est_maxwell_time}
\|\LU(T)\|_{\ee,\mm}
\leq
\|\LU_0\|_{\ee,\mm}
+
\int_0^T \|\LF(t)\|_{\ee,\mm} dt
\end{equation}
holds true. Owing to the regularity of $\LU$, simple manipulations then show that
we can rewrite the first line of \eqref{eq_ode_maxwell_time} as
\begin{equation}
\label{eq_var_maxwell_time}
(\dot \LU(t),v)_{\ee,\mm} + a(\LU(t),v) = (\LF(t),v)_{\ee,\mm} \quad \forall t \in [0,T]
\end{equation}
for all $v \in \LV(\Omega)$.

So far, we have defined solutions to \eqref{eq_maxwell_time_strong} in a
variational sense for sufficiently smooth initial data $\LU_0 \in \LVI$,
where the link between \eqref{eq_maxwell_weak} and \eqref{eq_var_maxwell_time}
is clear. This is not entirely sufficient since as previously explained, the
functional framework for the controllability method is set in $\LL(\Omega)$.
By density of $\LVI(\Omega)$ in $\LL(\Omega)$ however, estimate \eqref{eq_est_maxwell_time}
enables us to define, for any fixed $\LF$, the operator $\LU_0 \to \LU(T)$ for
all $\LU_0 \in \LL(\Omega)$ by continuity, thereby defining a continuous affine operator
mapping $\LL(\Omega)$ into itself. This observation is linked to the fact that
when $\LF \eq 0$, the operator $A$ is the infinitesimal generator of a $C_0$ semigroup
on $\LL(\Omega)$, see \cite{pazy_1983a}.

Although $\LU(T)$ can be defined for rough initial data
$\LU_0 \in \LL(\Omega)$, the corresponding solution $\LU$ only solves \eqref{eq_ode_maxwell_time}
in a very weak sense as we only have $\LU \in C^0(0,T;\LL(\Omega))$. In particular,
\eqref{eq_var_maxwell_time} does not hold. In the proofs below, we circumvent this difficulty
by establishing our results first for initial data in $\LVI(\Omega)$, and then extend them
to the general case by continuity owing to the dense inclusion $\LVI(\Omega) \subset \LL(\Omega)$.

Finally, we note that in view of \eqref{eq_adjoint}, for all
$\LU_0 \in \LVI(\Omega)$, there exists a unique
$\LU^\star \in C^1(0,T;\LL(\Omega)) \cap C^0(0,T,\LVI(\Omega))$ such that
\begin{equation}
\label{eq_adj_maxwell_time}
(v,\dot \LU^\star(t))_{\ee,\mm} + a(v,\LU^\star(t)) = 0 \quad \forall t \in [0,t]
\end{equation}
and $\LU^\star(0) = \LU_0$. Here, we can also extend
the notion of (weak) solutions to \eqref{eq_adj_maxwell_time} to any $\LU_0 \in \LL(\Omega)$,
as for \eqref{eq_var_maxwell_time}.

%% file: figures/boundary_condition_settings.tex
\begin{tikzpicture}[scale=2]

\draw (1.5,1.5) circle (1.5);

\draw (1,1) to[in=220,out=180] (1,2) to[out=180+220,in=90] (2,2) to[in=0,out=-90] (1,1);

\draw (1.5,1.5) node[anchor=west]{$\GP$};

\draw (2.5,1.) node{$\Omega$};

\draw (3,2) node[anchor=west]{$\GI$};

\end{tikzpicture}

%% file: 4-abstract-analysis.tex
\section{Properties of time-periodic solutions}
\label{section_perio}

Here, we introduce the key operators at involved in the controllability method.
We also discuss in detail the link between periodic solutions to time-dependent
Maxwell's equations \eqref{eq_maxwell_time_strong} and the time-harmonic solution
to \eqref{eq_maxwell_freq_strong}.

\subsection{Key operators}

First, we introduce the filtering and propagator operators,
which are the building blocks of the energy functional and
the associated CMCG method.

\subsubsection{Filtering}

Let $T \eq \omega/(2\pi)$ denote the period associated with the frequency $\omega$.
The filtering operator $F_\omega$ is defined by
\begin{equation}
\label{eq_filter}
F_\omega \LU \eq \frac{2}{T} \int_0^T \LU(t) e^{-i\omega t}dt
\end{equation}
for all $\LU \in C^0(0,T;\LL(\Omega))$. Clearly,
$F_\omega$ continuously maps $C^0(0,T;\LL(\Omega))$ into $L(\Omega)$.
and $C^0(0,T;\LV(\Omega))$ into $V(\Omega)$. In addition, when $\LU \in C^1(0,T;\LL(\Omega))$,
integration by parts easily shows that
\begin{equation}
\label{eq_filter_derivative}
F_\omega \dot \LU = i\omega F_\omega \LU + \frac{\omega}{\pi} \jmp{\LU}_T,
\end{equation}
where, for $\LW \in C^0(0,T,\LL(\Omega))$, we have introduced the notation
$\jmp{\LW}_T \eq \LW(T)-\LW(0)$.

\subsubsection{Propagators}
\label{section_perio_propagators}

Following the discussion in Section \ref{section_model_time},
if $\LU_0 \in \LVI(\Omega)$ and $\phi \in L(\Omega)$,
there exists a unique element $\LU \in C^1(0,T;\LL(\Omega)) \cap C^0(0,T;\LVI(\Omega))$
such that
\begin{equation}
\label{eq_def_propagator}
\left \{
\begin{array}{rcll}
(\dot \LU(t),v)_{\ee,\mm} + a(\LU(t),v)
&=&
(\Re (\phi \; e^{i\omega t}), v)_{\ee,\mm}
&
\forall v \in \LV, \quad t \in (0,T)
\\
\LU(0) &=& \LU_0, &
\end{array}
\right .
\end{equation}
and we define forward propagator $P_{\phi,\omega} \LU_0 \eq \LU(T)$.
When $\phi \eq 0$, we simply write $P_\omega \eq P_{0,\omega}$.

Similarly, we define a backward propagator. For $\LW_T \in \LVI(\Omega)$,
there exists a unique element $\LW \in C^1(0,T,\LL(\Omega)) \cap C^0(0,T,\LVI(\Omega))$
such that
\begin{equation}
\label{eq_def_backward_propagator}
\left \{
\begin{array}{rcll}
-(v,\dot \LW(t))_{\ee,\mm} + a(v,\LW(t))
&=&
0
&
\forall v \in \LV, \quad t \in (0,T)
\\
\LW(T) &=& \LW_T, &
\end{array}
\right .
\end{equation}
and we set $P_\omega^\star \LW_T \eq \LW(0)$.
Notice that $\LW$ is indeed well-defined, since the change of variable $\widetilde t \eq T-t$
transforms \eqref{eq_def_backward_propagator} into \eqref{eq_adj_maxwell_time}.
Together with \eqref{eq_adjoint}, this remark shows that the same time-stepping algorithm may be
used to compute $P_{\phi,\omega}$ and $P_\omega^\star$ simply by changing the sign of
the magnetic field.

Again, while the above definitions of $P_{\phi,\omega}$
and $P_\omega^\star$ require $\LVI(\Omega)$-regularity of the initial data,
semigroup theory allows us to extend the definitions of
$P_{\phi,\omega}$ and $P_\omega^\star$ as operators continuously mapping $\LL(\Omega)$
into itself \cite{pazy_1983a}.

Next, we remark that $P_\omega$ is linear, whereas $P_{\phi,\omega}$ is affine,
since
\begin{equation}
\label{eq_decomposition_propagator}
P_{\phi,\omega} \LU_0 = P_\omega \LU_0 + P_{\phi,\omega} 0
\qquad
\forall \LU_0 \in \LL(\Omega).
\end{equation}

\begin{lemma}
The operator $P_\omega^\star$ is the adjoint of $P_\omega$ for the
$\LL(\Omega)$ inner-product, i.e.
\begin{equation}
\label{eq_dual_propagators}
(P_\omega \LU_0,\LW_T)_{\ee,\mm}
=
(\LU_0,P_\omega^\star\LW_T)_{\ee,\mm}
\end{equation}
for all $\LU_0,\LW_T \in \LL(\Omega)$.
\end{lemma}

\begin{proof}
We only need to show \eqref{eq_dual_propagators} in $\LVI(\Omega)$;
the general case follows by density. Hence, we consider $\LU_0,\LW_T \in \LVI(\Omega)$
and denote by $\LU,\LW \in C^1(0,T,\LL(\Omega)) \cap C^0(0,T,\LVI(\Omega))$ the associated
solutions to \eqref{eq_def_propagator} and \eqref{eq_def_backward_propagator}. Owing to the
time-regularity of $\LU$ and $\LW$, integration by parts shows that
\begin{equation*}
\int_0^T (\dot \LU(t),\LW(t))_{\ee,\mm} dt
=
\left [
(\LU(t),\LW(t))_{\ee,\mm}
\right ]_0^T
-
\int_0^T (\LU(t),\dot \LW(t))_{\ee,\mm} dt,
\end{equation*}
which we rewrite as
\begin{equation}
\label{tmp_dual_propagator_0}
\int_0^T (\dot \LU(t),\LW(t))_{\ee,\mm} dt
+
\int_0^T (\LU(t),\dot \LW(t))_{\ee,\mm} dt
=
(P_\omega \LU_0,\LW_T)_{\ee,\mm}
-
(\LU_0,P_\omega^\star \LW_T)_{\ee,\mm}.
\end{equation}
The left-hand side of \eqref{tmp_dual_propagator_0} vanishes, since
\begin{multline*}
\int_0^T (\dot \LU(t),\LW(t))_{\ee,\mm} dt
+
\int_0^T (\LU(t),\dot \LW(t))_{\ee,\mm} dt
\\
=
\int_0^T (\dot \LU(t),\LW(t))_{\ee,\mm} + a(\LU(t),\LW(t)) dt
+
\int_0^T (\LU(t),\dot \LW(t))_{\ee,\mm} - a(\LU(t),\LW(t)) dt
\end{multline*}
which is zero due to \eqref{eq_def_propagator} and \eqref{eq_def_backward_propagator}.
\end{proof}

\subsubsection{Filtering of initial conditions}

If $\LU_0 \in \LL(\Omega)$ and $\phi \in L(\Omega)$,
we introduce the notation $F_{\phi,\omega} \LU_0 \eq F_\omega \LU$,
where $\LU \in C^0(0,T,\LL(\Omega))$ solves \eqref{eq_def_propagator}
in a weak sense, see \ref{section_model_time}. For
$\phi \eq 0$, we simply write $F_\omega \LU_0 \eq F_{0,\omega} \LU_0$.

\subsubsection{Energy functional}

Let $J: \LL(\Omega) \to \BbBR$ denote the ``energy functional''
\begin{equation}
\label{eq_energy_functional}
J(\LU_0) \eq \frac{1}{2} \|P_{\psi,\omega} \LU_0 - \LU_0\|_{\ee,\mm}^2
\quad
\forall \LU_0 \in \LL(\Omega).
\end{equation}
Using \eqref{eq_decomposition_propagator}, we can rewrite \eqref{eq_energy_functional} as
\begin{equation}
\label{eq_energy_functional_linear}
J(\LU_0) = \frac{1}{2} \|(I-P_\omega) \LU_0 - \LG\|_{\ee,\mm}^2
\quad
\forall \LU_0 \in \LL(\Omega),
\end{equation}
where $\LG \eq P_{\psi,\omega} 0$. Note that $J$ is continuous
over $\LL(\Omega)$ thanks to the discussions in Sections \ref{section_model_time}
and \ref{section_perio_propagators}.

\subsection{Structure of the minimizers}

For $U$, the (unique) time-harmonic solution to \eqref{eq_frequency_solution},
$\LU_0 \eq \Re U$ is a minimizer of $J$ since $J(\LU_0) = 0$. However, depending
on the boundary conditions, and properties of the right-hand sides, $\LU_0$
may not be the only minimizer of $J$. In this section, we
analyze the properties satisfied by the minimizers of $J$
and exhibit the structure of the minimization set. We also
identify situations in which the minimizer of $J$ is unique.

The starting point of our analysis is the following model decomposition result.

\begin{lemma}[Modal decomposition]
\label{lemma_modal_decomposition}
Let $\LU_0 \in \LVI(\Omega)$ satisfy $J(\LU_0) = 0$. Then, we have
\begin{equation}
\label{eq_fourier_U0}
\LU_0 = \Re \left (U_0 + U + \sum_{\ell \geq 2} U_\ell \right ),
\end{equation}
where $U_0 \in \ker a$, $U$ is the unique solution to \eqref{eq_frequency_solution},
and for $\ell \geq 2$, $U_\ell$ is an element of $V(\Omega)$ satisfying
\begin{equation}
\label{eq_higher_harmonics}
i\ell\omega (U_\ell,v) + a(U_\ell,v) = 0 \quad \forall v \in V(\Omega).
\end{equation}
\end{lemma}

\begin{proof}
Since the proof closely follows along the lines of \cite[Theorem 6]{tang_2020a},
we omit details for the sake of brevity. Consider $\LU_0 \in \LVI(\Omega)$ such
that $J(\LU_0) = 0$, and let $\LU \in C^1(0,T,\LL(\Omega)) \cap C^0(0,T,\LVI(\Omega))$
be the solution to \eqref{eq_def_propagator} with initial condition $\LU_0$ and
right-hand side $\psi$. By assumption, $J(\LU_0) = 0$ since $\LU$ is $T$-periodic.
Hence, we can expand $\LU$ in Fourier series as
\begin{equation}
\label{tmp_fourier_U0}
\LU(t) = \Re \left (\sum_{\ell \geq 0} U_\ell e^{i\ell\omega t} \right )
\quad
\forall t \in (0,T)
\end{equation}
where
\begin{equation}
\label{tmp_fourier_coef}
U_0 \eq \frac{1}{T} \int_0^T \LU(t) dt  \; \in V(\Omega),
\qquad
U_\ell \eq \frac{2}{T} \int_0^T \LU(t) e^{-i\ell\omega t} dt, \quad \ell \geq 1,
\end{equation}
Then, we obtain \eqref{eq_fourier_U0} by setting $t=0$ in \eqref{tmp_fourier_U0}.
After multiplying \eqref{eq_def_propagator} by $e^{-i\ell\omega t}$ and integrating
over $(0,T)$, we see that $U_0 \in \ker a$, $U_1 = U$, and that $U_\ell$ satisfies
\eqref{eq_higher_harmonics} for $\ell \geq 2$.
\end{proof}

Equipped with Lemma \ref{lemma_modal_decomposition}, we need a further understanding
of the kernel
\begin{equation*}
\ker a
\eq
\left \{
u \in V(\Omega) \; | \; a(u,v) = 0 \quad \forall v \in V(\Omega)
\right \}
\end{equation*}
and the space
\begin{equation*}
K(\Omega)
\eq
\left \{
[\be,\bh] \in V(\Omega)
\left |
\begin{array}{l}
\be \times \bn = \bh \times \bn = \zero \text{ on } \GI
\\
\curl \be = \curl \bh = \zero \text{ in } \Omega
\end{array}
\right .
\right \}
\end{equation*}
will play an important role. To characterize its structure, we introduce the set of gradients
$G(\Omega) \eq \grad H^1_\Gamma(\Omega,\mathbb C) \times \grad H^1_\GI(\Omega,\mathbb C)$
and its orthogonal complement (with respect to the $(\cdot,\cdot)_{\ee,\mm}$ inner-product)
$Z(\Omega) \eq G^\perp(\Omega)$, which consists of divergence-free functions. Then, we have
$K(\Omega) = G(\Omega) \oplus H(\Omega)$, where $H(\Omega) \eq K(\Omega) \cap Z(\Omega)$
is a ``cohomology'' space associated with $\Omega$. The structure of
$H(\Omega)$ is well-characterized \cite{fernandes_gilardi_1997a}. In particular,
it is finite-dimensional, and even trivial when $\Omega$ is simply-connected.
Similar properties hold for the real-valued counterparts of these spaces.

\begin{lemma}[Characterization of $\ker a$]
\label{lemma_kernel}
We have
\begin{equation*}
\ker a
=
\left \{
[\be,\bh] \in K(\Omega)
\; | \;
\be = \zero \text{ on } \Omega_{\sig}
\right \}.
\end{equation*}
\end{lemma}

\begin{proof}
Let $W \eq [\be,\bh] \in V(\Omega)$. For all smooth, compactly supported, vector
valued-function $\bphi \in \boldsymbol{\mathcal D}(\Omega)$, we have
\begin{equation*}
a([\be,\bh],[\bphi,0]) = (\sig \be,\bphi) + (\bh,\curl \bphi) = 0,
\qquad
a([\be,\bh],[0,\bphi]) = -(\be,\curl \bphi) = 0,
\end{equation*}
which implies that $\curl \bh = -\sig \be$ and $\curl \be = \zero$. As a consequence,
we have
\begin{align*}
0
&=
\Re a([\be,\bh],[\be,\bh])
\\
&=
(\sig \be,\be) + (\BY \be_\tau,\be_\tau)_{\GI} + (\BZ \bh_\tau,\bh\tau)_{\GI}
+
(\bh,\curl \be)-(\be,\curl \bh)
\\
&=
2(\sig \be,\be) + (\BY \be_\tau,\be_\tau)_{\GI} + (\BZ \bh_\tau,\bh\tau)_{\GI},
\end{align*}
from which we conclude that $\be \times \bn = \bh \times \bn = \zero$
on $\GI$ and $\be = \zero$ in $\Omega_{\sig}$. This last equality
also implies that $\curl \bh = \zero$.
\end{proof}

The first key result of this section applies to the case where the time-harmonic
problem is well-posed for all multiplies $\ell\omega$ of the original frequency
$\omega$.  It is an immediate consequence of Lemmas \ref{lemma_modal_decomposition}
and \ref{lemma_kernel} and of the decomposition of $K(\Omega)$
discussed above, so that its proof is omitted.

\begin{theorem}[Decomposition for well-posed problems]
\label{thm_decomposition_well_posed}
Assume that time-harmonic equations \eqref{eq_maxwell_weak} are well-posed for all
frequencies $\ell\omega$, $\ell \in \mathbb N^\star$. Then, we have
\begin{equation*}
\LU_0 = \Re \left ([\grad p,\grad q] + \theta + U\right )
\end{equation*}
where $p \in H^1_\Gamma(\Omega,\BbBC)$ and $q \in H^1_{\GI}(\Omega,\BbBC)$ and
$\theta \in H(\Omega)$.
\end{theorem}

Next, we show that if the right-hand side of the problem satisfies suitable
conditions, the ``stationary part'' $U_0$ of the minimizer must vanish.

\begin{theorem}[Decomposition of divergence-free minimizers]
\label{thm_decomposition_divergence_free}
Assume that $\psi \in K^\perp(\Omega)$ and that $\LU_0 \in \LV(\Omega) \cap \LK^\perp(\Omega)$.
Then, we have
\begin{equation*}
\LU_0 = \Re \left (U + \sum_{\ell \geq 2} U_\ell\right ).
\end{equation*}
\end{theorem}

\begin{proof}
Let $\LU$ be the time domain solution with initial condition $\LU_0$,
and introduce $[\BE_0,\BH_0] \eq \LU_0$ and $[\BE,\BH] \eq \LU$.
For any test functions $[\bv,\zero],[\zero,\bw] \in \LK(\Omega)$, we have
\begin{equation*}
(\ee\dot\BE,\bv)_{\widetilde \Omega_{\sig}} = (\mm\dot \BH,\bw)_\Omega = 0,
\end{equation*}
which implies that $[\BE(t),\BH(t)] \in \LK^\perp(\Omega)$. Therefore, $U_0 \in K^\perp(\Omega)$.
It follows that $U_0 \in K(\Omega) \cap K^\perp(\Omega)$ and hence, vanishes.
\end{proof}

We finally observe that if the assumptions of Theorems
\ref{thm_decomposition_well_posed} and \ref{thm_decomposition_divergence_free}
are both satisfied, we indeed have $\LU_0 = \Re U$. Since
$\LK^\perp(\Omega) = \LZ(\Omega) \cap \LH^\perp(\Omega)$, we see
that the assumptions on $\LU_0$ and $\psi$ in the statement of
\eqref{thm_decomposition_divergence_free} mean that these fields are
divergence-free and orthogonal to the (finite-dimensional) space $\LH(\Omega)$.
Note that this last requirement is null for simply connected domains, since
$\LH(\Omega) = \{0\}$ in this case. Similarly to \cite[Theorem 1]{grote_nataf_tang_tournier_2020a}
in the acoustic case, it is always possible to explicitly compute the time independent
components $[\grad p,\grad q]$ and $\theta$ by solving Poisson problems.

\subsection{Filtering of periodic solutions}
\label{section_filtering}

In the previous section, we exhibited the structure of the
minimizing set of $J$ using Fourier theory. As the filtering
operator essentially selects one specific Fourier mode, modal
decomposition \eqref{eq_fourier_U0} can be used to show how filtering acts on minimizers
of $J$. In fact, this technique was used in \cite{grote_nataf_tang_tournier_2020a} to show
that for any minimizer $\LU_0$ of $J$, we recover the time-harmonic solution $U$ after filtering.

Here, we develop an alternate proof technique, that actually
does not rely on the development of the previous section. This idea
appears to be new, and enables to quantify how well
initial conditions $\LU_0$ leading to ``approximately periodic''
time-dependent solution approximate the time-harmonic
solution $U$ after filtering. The proof improves similar concepts
used in \cite[Theorem 10]{tang_2020a} for the acoustic Helmholtz equation
formulated using a second-order in time framework.

\begin{theorem}[Alternate characterization of filtered solutions]
\label{theorem_filtering}
Let $\phi \in L(\Omega)$. Then, for all $\LU_0 \in \LL(\Omega)$,
we can characterize $F_\omega \LU_0$ as the unique element of $V(\Omega)$ such that
\begin{equation}
\label{eq_characterization_filtering}
i\omega(F_{\phi,\omega} \LU_0,v)_{\ee,\mm} + a(F_{\phi,\omega} \LU_0,v)
=
(\phi,v)_{\ee,\mm} + \frac{\omega}{\pi} (\LU_0-P_\omega\LU_0,v)_{\ee,\mm}
\end{equation}
for all $v \in V(\Omega)$. As a direct consequence, we have
\begin{equation}
\label{eq_mismatch_stability}
\enorm{U-F_{\psi,\omega} \LU_0}
\leq
\frac{\omega}{\pi} \Cstab \|(I-P_{\psi,\omega})\LU_0\|_{\ee,\mm}.
\end{equation}
for all $\LU_0 \in \LV(\Omega)$.
\end{theorem}

\begin{proof}
We first discuss the case where $\LU_0 \in \LVI(\Omega)$. Thus,
let $\LU$ be as in \eqref{eq_def_propagator} with initial condition
$\LU_0$ and right-hand side $\phi \in L(\Omega)$.
For all $v \in \LV(\Omega)$, we have
\begin{equation}
\label{tmp_proof_mismatch}
\frac{2}{T}
\int_0^T \left \{(\dot \LU,v)_{\ee,\mm} + a(\LU,v)
\right \} e^{-i\omega t} dt
=
\frac{2}{T} \int_0^T (\Re (\phi e^{i\omega t}), v)_{\ee,\mm} e^{-i\omega t} dt.
\end{equation}
Since $\ee,\sig,\mm$ and $v$ are time-independent, we can write
\begin{equation*}
\frac{2}{T}
\int_0^T \left \{(\dot \LU,v)_{\ee,\mm} + a(\LU,v)
\right \} e^{-i\omega t} dt
=
(F_\omega \dot \LU,v)_{\ee,\mm} + a (F_\omega \LU,v),
\end{equation*}
and \eqref{eq_filter_derivative} shows that
\begin{equation*}
\frac{2}{T}
\int_0^T \left \{(\dot \LU,v)_{\ee,\mm} + a(\LU,v)
\right \} e^{-i\omega t} dt
=
i\omega (F_\omega \LU,v)_{\ee,\mm} + a(F_\omega \LU,v)
+
\frac{\omega}{\pi} (\jmp{\LU}_T,v)_{\ee,\mm}.
\end{equation*}
Similarly, since $\phi$ is time-independent, we have
\begin{equation*}
\frac{2}{T} \int_0^T (\Re (\phi \; e^{i\omega t}), v)_{\ee,\mm} e^{-i\omega t} dt
=
(\phi,v)_{\ee,\mm},
\end{equation*}
and as a result
\begin{equation*}
i\omega (F_\omega \LU,v)_{\ee,\mm} + a(F_\omega \LU,v)
=
(\phi,v)_{\ee,\mm} - \frac{\omega}{\pi} (\jmp{\LU}_T,v)_{\ee,\mm},
\end{equation*}
so that \eqref{eq_characterization_filtering} follows whenever $\LU_0 \in \LVI(\Omega)$,
recalling that $F_{\phi,\omega} \LU_0 \eq F_\omega \LU$
and $\jmp{\LU}_T \eq P_{\phi,\omega} \LU_0 - \LU_0$.

For the general case where $\LU_0 \in \LL(\Omega)$, we first observe that we may
equivalently rewrite \eqref{eq_characterization_filtering} as
\begin{equation}
\label{tmp_characterization_L2}
F_{\phi,\omega} \LU_0
=
S_\omega\left (\phi + \frac{\omega}{\pi}(I-P_\omega) \LU_0 \right ).
\end{equation}
At that point, identity \eqref{tmp_characterization_L2} is already established in $\LVI(\Omega)$.
But then, since \eqref{tmp_characterization_L2} involves continuous operators
from $L(\Omega)$ into itself, the density of $\LVI(\Omega)$ into $L(\Omega)$ implies
the general case.

To conclude the proof, letting $\phi = \psi$ and recalling the definition
\eqref{eq_frequency_solution} of $U$, we obtain
\begin{equation*}
i\omega(U-F_{\psi,\omega} \LU_0,v)_{\ee,\mm}
+
a(U-F_{\psi,\omega} \LU_0,v)
=\frac{\omega}{\pi}((P_{\psi,\omega}-I)\LU_0,v)_{\ee,\mm},
\end{equation*}
so that \eqref{eq_mismatch_stability} follows from \eqref{eq_frequency_stability}.
\end{proof}

Using \eqref{eq_maxwell_weak}, we may rewrite \eqref{eq_characterization_filtering}
in compact form as
\begin{equation}
\label{eq_characterization_filtering_compact}
F_\omega \LU_0 = \frac{\omega}{\pi} S_\omega \circ (I-P_\omega) \LU_0
\qquad \forall \LU_0 \in \LL(\Omega).
\end{equation}
Taking again advantage of the similarity between the
original and adjoint problems, we can also show that
\begin{equation}
\label{eq_characterization_filtering_compact_adjoint}
F_\omega \LW_T = \frac{\omega}{\pi} S_\omega^\star \circ (I-P_\omega^\star) \LW_T
\qquad
\forall \LW_T \in \LL(\Omega).
\end{equation}

Stability estimate \eqref{eq_mismatch_stability} is of particular interest,
since it shows that filtering ``nearly periodic'' solutions yields
good approximations of the time-harmonic solution. It also suggests that the
misfit $\LU_0-P_{\psi,\omega} \LU_0$ may be used as a stopping criterion for iterative
methods, but the dependency on the frequency must be taken into account.

%% file: 5-controllability-method.tex
\section{Controllability Method}
\label{section_CMCG}

In this section, we build upon the results of the previous section
to introduce our controllability method, that we couple with a conjugate
gradient minimization algorithm.

We seek an initial condition $\LU_0 \in \LL(\Omega)$ satisfying $P_{\psi,\omega} \LU_0 = \LU_0$,
or maybe more explicitly, such that
\begin{equation}
\label{eq_controlability}
(I-P_\omega)\LU_0 = \LG,
\end{equation}
where $P_{\psi,\omega}$, $P_{\psi}$ and $\LG$ are respectively introduced at
\eqref{eq_def_propagator}, \eqref{eq_decomposition_propagator}
and \eqref{eq_energy_functional_linear}.
Clearly, $\LU_0 \eq \Re U$ is one solution to \eqref{eq_controlability}
but it may not be unique. Nevertheless, we always have $U = F_{\psi,\omega} \LU_0$.
In addition, estimate \eqref{eq_mismatch_stability}
implies that for any approximate solution $\LU_0$ to \eqref{eq_controlability},
$F_\omega \LU_0$ is an approximate solution to \eqref{eq_frequency_solution}.

\subsection{The conjugate gradient method}

After space discretization, \eqref{eq_controlability} corresponds to a finite-dimensional
linear system. In principle, the matrix corresponding to $P_\omega$ could therefore be
(approximately) assembled by running a time-domain solver for one period for every possible initial
conditions.  However, this approach is prohibitively expensive in practice. Instead, we opt
for the matrix-free conjugate gradient iteration, which only requires
evaluating $P_\omega \LU_0$ for a limited number of initial conditions.

We thus reformulate controllability equation \eqref{eq_controlability} as the optimization problem
\begin{equation}
\label{eq_controllability_formulation}
\min_{\LU_0 \in \LL(\Omega)} J(\LU_0),
\end{equation}
where $J$ is the energy functional introduced in \eqref{eq_energy_functional}.
From \eqref{eq_energy_functional_linear}, we recall that
$J$ corresponds to a ``standard'' quadratric form and, as result,
its gradient and Hessian are easily derived. The proof of the result below is omitted,
as it follows from standard algebraic manipulations.

\begin{theorem}[Structure of the energy functional]
For all $\LU_0,\LV_0 \in \LL(\Omega)$, we have
\begin{align*}
J(\LU_0+\LV_0)
=
J(\LU_0)
&+
\Re ((I-P_\omega^\star)(I-P_\omega)\LU_0-(I-P_\omega^\star)\LG,\LV_0)_{\ee,\mm}
\\
&+
\frac{1}{2}
((I-P_\omega)\LV_0,(I-P_\omega)\LV_0)_{\ee,\mm}.
\end{align*}
It follows that
\begin{equation}
\label{eq_functional_gradient}
J'(\LU_0) =  (I-P_\omega^\star)(I-P_\omega) \LU_0 - (I-P_\omega^\star) \LG
\end{equation}
and
\begin{equation}
\label{eq_functional_hessian}
\left (J''(\LU_0)\right )(\LV_0,\LV_0)
=
\|(I-P_\omega) \LV_0\|_{\ee,\mm}^2.
\end{equation}
\end{theorem}

Next, we show that $J$ is continuous, uniformly Lipschitz, and strongly convex
over the quotient space $\LL(\Omega) / \ker F_\omega$. These properties ensure
the uniqueness of the minimizer of $J$ up to an element of $\ker F_\omega$ and
also implies the convergence of gradient-based algorithms \cite{ciarlet_1989a}.

\begin{theorem}[Convexity of energy functional]
For $\LU_0 \in \LL(\Omega)$, we have
\begin{equation}
\label{eq_functional_continuous}
J(\LU_0)
=
\frac{1}{2}
\left \|\frac{\pi}{\omega} S_\omega^{-1}F_\omega \LU_0 - \LG \right \|_{\ee,\mm}^2.
\end{equation}
In addition, for all $\LU_0,\LV_0 \in \LL(\Omega)$, the estimates
\begin{equation}
\label{eq_functional_lipschitz}
\|J'(\LU_0)-J'(\LV_0)\|_{\ee,\mm}
\leq
\frac{\omega^2}{\pi^2} \enorm{F_\omega(\LU_0-\LV_0)}
\end{equation}
and
\begin{equation}
\label{eq_functional_convex}
\left (J''(\LU_0)\right )(\LV_0,\LV_0)
\geq
\frac{\pi^2}{\omega^2}\frac{1}{\Cstab^2} \enorm{F_\omega \LV_0}^2
\end{equation}
hold true.
\end{theorem}

\begin{proof}
Identity \eqref{eq_functional_continuous} is a direct
consequence of \eqref{eq_characterization_filtering}.
Then, estimate \eqref{eq_functional_lipschitz} follows
from \eqref{eq_functional_gradient},
characterizations \eqref{eq_characterization_filtering_compact}
and \eqref{eq_characterization_filtering_compact_adjoint} of
$(I-P_\omega)$ and $(I-P_\omega^\star)$, and the continuity estimate
\eqref{eq_frequency_continuity}. Finally, we obtain convexity estimate
\eqref{eq_functional_convex} from \eqref{eq_functional_hessian},
\eqref{eq_characterization_filtering} and \eqref{eq_frequency_stability}.
\end{proof}

This result is to be compared with \cite[Theorem 3]{bardos_rauch_1994a},
where a convexity result is established under specific assumptions on the
spectrum. The use of the filtering allows to bypass this limitation.

In practice, it is not necessary to introduce the quotient space
$\LL(\Omega) / \ker F_\omega$. Indeed, a careful examination of standard convergence
proofs (see, e.g., \cite[Theorem 8.4.4]{ciarlet_1989a}) shows that properties
\eqref{eq_functional_lipschitz} and \eqref{eq_functional_convex} are sufficient to ensure
the convergence of $F_{\psi,\omega} \LU_0^{(\cgit)}$ to $U$ starting from any initial guess
$\LU_0^{(0)} \in \LL(\Omega)$, where $\LU_0^{(\cgit)}$ denotes a minimizing sequence.
In addition, a reduction factor of the form
\begin{equation*}
\enorm{U-F_{\psi,\omega} \LU_0^{(\cgit+1)}}
\leq
\left (1-\Cstab^{-4}\right ) \enorm{U-F_{\psi,\omega} \LU_0^{(\cgit)}}
\end{equation*}
can be obtained.

Among the possible gradient descent techniques, we select the usual CG iteration
(see \cite[Section 8.5]{ciarlet_1989a}) to solve \eqref{eq_controllability_formulation}.

\subsection{Discretization}

In our computations, we use an upwind-flux discontinuous Galerkin
method to discretize Maxwell's equations \eqref{eq_def_propagator} and
\eqref{eq_def_backward_propagator} in space, while explicit Runge-Kutta schemes are
employed for time integration. We restrict our numerical experiments to
two-dimensional examples, and the required notation is briefly presented below.

\subsubsection{Two-dimensional setting}

Here, we consider two-dimensional Maxwell's equations in a bounded domain
$\Omega\subset\BbBR^2$. Specifically, we consider three-dimensional Maxwell's equations
\eqref{eq_maxwell_freq_strong} in the domain $\Omega \times I$ for some interval $I$,
under the assumption that the electromagnetic field $(\be,\bh)$ does not depend on the
third space variable. There are two uncoupled polarizations, and we focus on the
``transverse magnetic'' case where $\bh = (\bh_1,\bh_2,0)$ and $\be = (0,0,\be_3)$.
The other polarization can be dealt with similarly by swapping the roles of $\bh$ and $\be$.
Employing the notation $\bh$ for the 2D vector gathering the magnetic field component and $e$
for the only non-zero component of the electric field. This, time-harmonic Maxwell's equations
reduce to
\begin{equation}
\label{eq_Maxwell_TM}
 \left\{
\begin{array}{rclll}
i\omega\varepsilon e + \sigma e + \scurl \bh &=& j & \text{ in } \Omega,
\\
i\omega \mu \bh - \vcurl e &=& \zero & \text{ in } \Omega,
\\
e &=& 0 & \text{ on } \GP,
\\
e + Z \bh_\tau &=& g &\text{ on } \GI,
\end{array}
\right.
\end{equation}
where $\varepsilon,\sigma,\mu$ and $Z$ are now scalar-valued functions,
and the two-dimensional curl operators are given by
\begin{equation*}
\scurl \bv = \partial_1 \bv_2 - \partial_2 \bv_1
\qquad
\vcurl v = (\partial_2 v,-\partial_1 v)
\end{equation*}
for any vector-valued and scalar-valued function $\bv$ and $v$.

The corresponding time-domain Maxwell's equations are given by
\begin{subequations}
\label{eq_time_maxwell_2D}
\begin{equation}
\label{eq_time_maxwell_2D_volume}
\left \{
\begin{array}{rcl}
\varepsilon \dot E + \sigma E + \scurl \BH &=& J,
\\
\mu \dot \BH - \vcurl E &=& \zero,
\end{array}
\right.
\end{equation}
in $\Omega$ and
\begin{equation}
\label{eq_time__maxwell_2D_surface}
\left \{
\begin{array}{rcll}
E &=& 0 & \text{ on } \GP,
\\
E + Z \BH \times \bn &=& G & \text{ on } \GI,
\end{array}
\right.
\end{equation}
\end{subequations}
for all $t \in [0,T]$.

%
%

\subsubsection{Discontinuous Galerkin discretization}

Following \cite{fezoui_lanteri_lohrengel_piperno_2005a,hesthaven_warburton_2002a},
we discretize \eqref{eq_time_maxwell_2D} with a first-order discontinuous Galerkin
(DG) method. The computational domain $\Omega$ is thus partitioned into a mesh
$\CT_h$ consisting of triangular elements $K$. For any element $K \in \CT_h$,
$\rho_K$ denote the diameter of the largest circle contained in $K$.

For the sake of simplicity, we assume that $\CT_h$ is conforming in the sense that
the intersection $\overline{K_+} \cap \overline{K_-}$ of two distinct elements
$K_\pm \in \CT_h$ is either empty, a single vertex, or a full face of both elements.
Note that the considered DG method is very flexible, and can, in principle, accommodate
non-conforming meshes with hanging nodes and/or different types of elements.

Next, we denote by $\CF_h$ the set of faces associated to $\CT_h$,
and we assume that each boundary face $F \in \CF_h$ with $F \subset \partial \Omega$
either entirely belongs to $\GI$ or $\GP$. The sets $\CFIh,\CFPh \subset \CF_h$
gather those faces respectively lying in $\GI$ and $\CP$, whereas $\CFih$ gathers
the remaining ``interior'' faces. We associate with each face $F \in \CF_h$ a fixed
normal unit normal vector $\bn_F$ chosen such that $\bn_F = \bn$ when $F \subset \partial \Omega$.
For internal faces, the orientation is arbitrary. We also employ the
notation $\bt_F$ for the unit tangential to $F$ obtained from $\bn_F$
by a $+\pi/2$ rotation.

For a given integer $q \in \mathbb N$, $\CP_q(\CT_h)$
stands for scalar-valued functions $v: \Omega \to \mathbb R$
such that $v|_K$ is a polynomial of degree less than or equal to $q$
for all $K \in \CT_h$. Note that the elements of $\CP_q(\CT_h)$ are,
in general, discontinuous across the faces $F \in \CF_h$ of the mesh.
Similarly $\BCP_q(\CT_h)$ is the space of vector-valued functions
$\bv \eq (\bv_1,\bv_2): \Omega \to \mathbb R^2$ such that $\bv_1,\bv_2 \in \CP_q(\CT_h)$.

If $v \in \CP_q(\CT_h)$ and $F \in \CFih$, the notations
\begin{equation*}
\avg{v}_F \eq v_+|_F + v_-|_F
\qquad
\jmp{v}_F \eq v_+|_F (\bn_+ \cdot \bn_F) + v_-|_F (\bn_- \cdot \bn_F)
\end{equation*}
stand for the usual average and jump of $v$ across $F$, where we used
$v_\pm \eq v|_{K_\pm}$ and $\bn_\pm = \bn_{K_\pm}$,
for any to elements $K_-$ and $K_+$ of $\CT_h$ such that $F = \partial K_- \cap \partial K_+$.
For external faces, we simply set $\avg{v}_F \eq \jmp{v}_F \eq v|_F$.
In addition, if $\bw \in \BCP_q(\CT_h)$ the same notations have to be understood
component-wise.

Given $E_{h,0} \in \CP_q(\CT_h)$ and $\BH_{h,0} \in \BCP_q(\CT_h)$,
the semi-discrete DG scheme consists in finding $E_h(t) \in \CP_q(\CT_h)$ and
$\BH_h(t) \in \BCP_q(\CT_h)$ by solving the system of ODE
for $t \in (0,T)$,
\begin{equation}
\label{eq_DG_semidiscrete}
\left \{
\begin{array}{rcll}
(\eps \dot E_h(t),v_h)_\Omega
+
(\sigma E_h(t),v_h)_\Omega
+
(\BH_h(t),\vcurl v_h)_\Omega
+
(\widehat \BH_h(t) \times \bn_F,\jmp{v_h})_{\CF_h}
&=&
(J(t),v_h)
\\
(\mm \dot \BH_h(t),\bw_h)_\Omega
+
(E_h(t),\scurl \bw_h)_\Omega
+
(\widehat E_h(t),\jmp{\bw_h} \times \bn_F)_{\CF_h}
&=&
0
\end{array}
\right .
\end{equation}
for all $v_h \in \CP_q(\CT_h)$ and $\bw_h \in \CP_q(\CT_h)$,
with initial conditions $E_h(0) = E_{h,0}$ and $\BH_h(t) = \BH_{h,0}$.
In \eqref{eq_DG_semidiscrete}, $(\cdot,\cdot)_{\CF_h} \eq \sum_{F \in \CF_h} (\cdot,\cdot)_F$,
while $\widehat E_h(t)$ and $\widehat \BH_h(t)$
are the upwind ``numerical fluxes''
\begin{equation*}
\widehat E_h|_F
\eq
\frac{1}{\avg{Y_{\rm flux}}} \left (
\avg{Y_{\rm flux} E_h}_F + \frac{1}{2} \jmp{\BH_h}_F \times \bn_F
\right )
\quad
\widehat
\BH_h|_F
\eq
\frac{1}{\avg{Z_{\rm flux}}} \left (
\avg{Z_{\rm flux}\BH_h}_F - \frac{1}{2} \jmp{E_h}_F \bt_F
\right ),
\end{equation*}
where $Z_{\rm flux} \eq \sqrt{\mu/\varepsilon}$, $Y_{\rm flux} = 1/Z_{\rm flux}$,
whenever $F \in \CFih$. For the remaining faces, we set
\begin{equation*}
\widehat E_h|_F
\eq
0
\qquad
\widehat
\BH_h|_F
\eq
-Y E_h \bt_F + \BH_h
\end{equation*}
when $F \in \CFPh$ and
\begin{equation*}
\widehat E_h|_F
\eq
\frac{1}{2} \left (E_h + Z\BH_h \times \bn + G \right )
\qquad
\widehat
\BH_h|_F
\eq
\frac{Y}{2} \left ( Z\BH_h - E_h \bt_F - G\bt_F \right )
\end{equation*}
if $F \in \CFIh$. This choice introduces some numerical dissipation,
leading to stable discretizations when coupled with Runge-Kutta time-integration.

To simplify further discussions, we introduce the compact notation
$\LU_h(t) \eq (E_h(t),\BH_h(t))$, and we
denote by $\VUh(t)$ the coefficients of $\LU_h(t)$ expanded in the
nodal basis of $\CP_q(\CT_h)$, to rewrite \eqref{eq_DG_semidiscrete} as
\begin{equation*}
\mathbb M \dot {\mathbb U}_h(t) + \mathbb K \VUh(t)
=
\Re \left (\mathbb M \mathbb J e^{i\omega t}\right ),
\end{equation*}
where $\mathbb M$ and $\mathbb K$ are the usual mass and stiffness matrices.
A key asset of DG discretizations is that $\mathbb M$ is block-diagonal, so that
the inverting $\mathbb M^{-1}$ is cheap. Hence, we may reformulate
the above ODE system as
\begin{equation}
\label{eq_ode_system}
\dot {\mathbb U}_h(t)
=
\Phi(t,\VUh(t)),
\qquad
\Phi(t,\VUh(t))
\eq
\Re \left (\mathbb J e^{i\omega t}\right )
+
\mathbb{B} \VUh(t),
\qquad
\mathbb B \eq \mathbb M^{-1} \mathbb K.
\end{equation}

\subsection{Time integration scheme}

We integrate \eqref{eq_ode_system} using a standard second-order explicit Runge-Kutta (RK2)
method with $\CP_1$ elements, or a fourth-order explicit Runge-Kutta (RK4) method with $\CP_3$
elements. Both are stable under a ``CFL condition'' on the time-step $\delta t$:
\begin{equation}
\label{eq_cfl}
\delta t \leq c_q \min_{K \in \CT_h} \left (\sqrt{\mu_K\varepsilon_K} \rho_K\right ),
\end{equation}
where the constant $c_q$ only depends on the polynomial degree $q$ and the shape-regularity
of the mesh. In our computations, we use $c_1 \eq 0.24$ and $c_3 \eq 0.12$, which we empirically
found to be near the stability limit.

We thus select a time-step $\delta t \eq T/M$, where $M$ is the smallest positive integer
such that \eqref{eq_cfl} holds, and iteratively compute approximation $\LU_{h,m}$ to $\LU_h(t_m)$
for $1 \leq m \leq M$, where $t_m \eq m \delta t$. Since there are no ``physical''
initial conditions, we are free to choose the initial condition as
piecewise polynomial function and therefore, there are no requirements to interpolate
or project the initial condition to define $\LU_{h,0}$ and the associated dof
vector $\mathbb U_{h,0}$. We either use the RK2 or the RK4 scheme to
compute $\mathbb U_{h,m+1}$ from $\mathbb U_{h,m}$. Both time integration schemes are
standard but for the sake of completeness, there are briefly listed in Algorithms \ref{al_rk2}
and \ref{al_rk4}.

\begin{algorithm}[t]
\caption{Explicit second-order Runge-Kutta (RK2) method}
\label{al_rk2}
\begin{algorithmic}[1]
\REQUIRE{$\mathbb U_{h,m}$ an approximation of $\VUh(t_m)$, $m\ge0$}
\STATE ${\mathbb K}_{h,1}  \ \eq\  \Phi(t_m,\mathbb U_{h,m})$
\STATE ${\mathbb K}_{h,2}  \ \eq\  \Phi(t_m+(\delta t/2),\mathbb U_{h,m} + (\delta t/2) {\mathbb K}_{h,1})$
\RETURN $\mathbb U_{h,m+1} \ \eq\  \mathbb U_{h,m} + \delta t {\mathbb K}_{h,2}$
\end{algorithmic}
\end{algorithm}

\begin{algorithm}[t]
\caption{Explicit fourth-order Runge-Kutta (RK4) method}
\label{al_rk4}
\begin{algorithmic}[1]
\REQUIRE{$\mathbb U_{h,m}$ an approximation of $\VUh(t_m)$, $m\ge0$}
\STATE ${\mathbb K}_{h,1}  \ \eq \  \Phi(t_m,\mathbb U_{h,m})$
\STATE ${\mathbb K}_{h,2}  \ \eq \  \Phi(t_m+(\delta t/2),\mathbb U_{h,m} + (\delta t/2) {\mathbb K}_{h,1})$
\STATE ${\mathbb K}_{h,3}  \ \eq \  \Phi(t_m+(\delta t/2),\mathbb U_{h,m} + (\delta t/2) {\mathbb K}_{h,2})$
\STATE ${\mathbb K}_{h,4}  \ \eq \  \Phi(t_m+\delta t,\mathbb U_{h,m} + \delta t {\mathbb K}_{h,3})$
\RETURN $\mathbb U_{h,m+1} \ \eq \ \mathbb U_{h,m}+(\delta t/6)\big({\mathbb K}_{h,1}+2{\mathbb K}_{h,2}+2{\mathbb K}_{h,3}+{\mathbb K}_{h,4}\big)$
\end{algorithmic}
\end{algorithm}

\subsection{Implementation of the filtering}

In this section, we briefly discuss the implementation of the filtering operator $F_\omega$
defined in \eqref{eq_filter}. For the RK2 scheme, we may simply employ the trapezoidal rule
\begin{equation}
\label{eq_filtering_RK2}
F_\omega \mathbb U_h
\simeq
\frac{\delta t}{2} \sum_{m=1}^M \left (
\mathbb U_{h,m-1} e^{-i\omega t_{m-1}} + \mathbb U_{h,m} e^{-i\omega t_m}
\right ),
\end{equation}
since it is second-order accurate. The situation is slightly more delicate for the RK4
scheme, as employing \eqref{eq_filtering_RK2} would deteriorate the convergence rate
of the method. Instead, we employ a method based on Hermite interpolation. This method
is especially efficient, because the RK algorithm computes the vectors $\Phi(t,\mathbb U_{h,m})$
anyways which are natural approximations to $\dot {\mathbb U}_{h,m}$. We thus let
\begin{equation*}
\mathbb I_{h,m}(t)
\eq
             \mathbb U_{h,{m-1}}  p_{00}(t) + \mathbb U_{h,m}           p_{01}(t)
+
\Phi(t_{m-1},\mathbb U_{h,{m-1}}) p_{10}(t) + \Phi(t_m,\mathbb U_{h,m}) p_{11}(t),
\end{equation*}
where the Hermite polynomials $p_{ij}$ are the only elements of $\CP_3(t_{m-1},t_m)$
satisfying $p_{ij}^{(\ell)}(t_{m+k}) = \delta_{ik} \delta_{j\ell}$ for $0 \leq k,\ell \leq 1$.
Since Hermite polynomials are explicitly available, we can evaluate
\begin{equation*}
\xi_{ij} \eq \int_{t_{m-1}}^{t_m} p_{ij}(t) e^{-i\omega t} dt
\end{equation*}
analytically, which yields
\begin{align}
\label{eq_filtering_RK4}
F_\omega \mathbb U_h
&\simeq
\sum_{m=1}^M \int_{t_{m-1}}^{t_m} \mathbb I_{h,m}(t) e^{-i\omega t}
\\
\nonumber
&=
             \mathbb U_{h,{m-1}}  \xi_{00} + \mathbb U_{h,m}           \xi_{01}
+
\Phi(t_{m-1},\mathbb U_{h,{m-1}}) \xi_{10} + \Phi(t_m,\mathbb U_{h,m}) \xi_{11}.
\end{align}

We emphasize that \eqref{eq_filtering_RK2} and \eqref{eq_filtering_RK4}
only require the solutions $\mathbb U_{h,{m-1}}$ and $\mathbb U_{h,m}$.
In fact, we can easily reformulate the above formula to only require
$\mathbb U_{h,m}$ at a single time, and this readily compute $F_\omega \mathbb U_h$
on the fly.

%% file: 6-numerical-examples.tex
\section{Numerical examples}
\label{section_numer}

This section gathers numerical examples where we compare our CMCG algorithm
against a limiting amplitude principle, where ``naive'' time-stepping is employed
until convergence. The latter algorithm is denoted by FW (for full wave).
We utilize the DG method described in Section \ref{section_CMCG} in both cases,
so that a fair measure of the cost is the number of periods that need to be simulated
to reach a given accuracy level. We chose to start both algorithm with $\LU_0^{(0)} = 0$
in all the considered experiments. It is known that this strategy is not optimal, since
transient right-hand sides generally improves the perfomance of FW, and the convergence of CMCG
can be accelerated, if it is applied after a ``run-up'' phase of a few FW iterations
(see, e.g.  \cite{bristeau_glowinski_periaux_1994a,tang_2020a}). Nevertheless,
we restrict ourselves to zero initialization for a fair comparison.

Another question we address is the comparison of the solution obtained after convergence
of the CMCG or FW method against the solution given by the same frequency-domain
DG scheme. In this case we solve the linear system
$(i\omega \mathbb M + \mathbb K) \mathbb U_h = \mathbb M \mathbb J_h$,
with the direct solver implemented in the software package {\tt MUMPS}
\cite{amestoy_ashcraft_boiteau_buttari_lexcellent_weisbecker_2015a,amestoy_duff_lexcellent_2000a}.
We use the notation FS (frequency solver) to refer to this solution.
This is a subtle point, because the CMCG and FW algorithm will converge to
a (slightly) different approximation, due to the error from time discretization.


Whenever the exact solution is available, we chose the mesh $\CT_h$ and
polynomial degree $q$ so that the FS relative error, measured as
\begin{equation*}
\text{error} \eq \|U-U_h\|_{\ee,\mm}/\|U\|_{\ee,\mm},
\end{equation*}
where $U$ is the exact solution and $U_h$ the FS solution,
is of the order of a few percents, which seems realistic for typical applications.
For the CMCG and FW method, the main figure of merit is then the relative error
\begin{equation*}
\text{error} \eq \|U - F_{\psi,\omega} \LU_{0,h}^{(\cgit)}\|_{\ee,\mm}/\|U\|_{\ee,\mm},
\end{equation*}
where $\LU_{h,0}^{(\cgit)}$ is the current iterate in the CMCG or FW algorithm. Specifically
$\LU_{h,0}^{(\cgit)}$ denotes the solution obtained after $\cgit$ iterations of the CMCG algorithm,
or the solution in the FW algorithm after simulating $\cgit$ periods. Note that CMCG requires
twice many time-periods to compute $\LU_{h,0}^{(\cgit)}$ as FW, which is accounted for in
the graphs below. In the last experiment, where the analytical solution is not available,
we monitor
\begin{equation*}
\text{error} \eq \|U_h - F_{\psi,\omega} \LU_{0,h}^{(\cgit)}\|_{\ee,\mm}/\|U_h\|_{\ee,\mm},
\end{equation*}
when comparing CMCG against FW.

In all examples we set $\sigma \eq 0$, $\mu \eq 1$, and $Z \eq 1$.
For $\theta \in [0,2\pi)$, we denote by $\bd_\theta \eq (\cos \theta,\sin \theta)$
the direction associated with $\theta$ and $\xi_\theta(\bx) \eq e^{i\omega \bd \cdot \bx}$
($\bx \in \BbBR^2$) is the plane wave travelling along the direction $\bd$.

Sometimes, we employ structured meshes based on Cartesian grids.
In this case, an ``$N \times M$ Cartesian mesh'' is obtained by starting
from a grid of $N \times M$ rectangles and then dividing each rectangle
into four triangles by joining each of its vertices with its barycentre.

\subsection{Plane wave in free space}

In this experiment, we set $\theta = 45^o$ and consider the propagation of a plane wave,
traveling along the direction $\bd_\theta$ in the square $\Omega \eq (0,1)^2$.
A Silver-M\"uller absorbing boundary condition is imposed on the whole boundary,
so that $\GI \eq \partial \Omega$ and $\GP \eq \emptyset$. We set $\eps \eq 1$,
$j \eq 0$ and $g = \grad \xi_\theta \cdot \bn + i\omega \xi_\theta$.
The solution then reads $(e,\magn) \eq (\xi_\theta, \xi_\theta \bd^\perp)$, with
$\bd^\perp \eq (-\sin \theta,\cos \theta)$. 

We consider the two frequencies $\omega = 10\pi$ and $40\pi$.
We employ a $32 \times 32$ Cartesian meshes in both cases with
$\CP_1$ elements for $\omega=10\pi$, and $\CP_3$ elements for $\omega=40\pi$.
Figure \ref{figure_planewave_error} shows the evolution of the error.
In this particular experiment, FW outperforms CMCG. When using $\CP_1$ elements, the error
achieved by both FW and CMCG is indistinguishable from the FS error. On the other hand,
the error slightly increases in both FW and CMCG when using $\CP_3$ elements.

\input{figures/planewave_error}

\subsection{Half open waveguide}

We now consider a rectangular domain $\Omega \eq (0,4) \times (0,1)$,
where the bottom, top and left sides are perfectly conducting,
while an impedance boundary condition is imposed on right side. Hence,
we have $\GP \eq (0,4) \times \{0,1\} \cup \{0\} \times (0,1)$ and
$\GI \eq \{4\} \times (0,1)$. Then, we solve \eqref{eq_Maxwell_TM} with
$\varepsilon \eq 1$, $j \eq 0$, $g \eq \xi_\theta$ and $\theta = 30^o$.

We obtain a semi-analytical solution by first performing the expansion
\begin{equation}
\label{eq_series_waveguide}
e = \sum_{n \geq 0} e_n(\bx_1) \sin(n\pi \bx_2),
\end{equation}
that is justified by the fact that the top and bottom boundary conditions
are ``Dirichlet-like''. Then, $e_n$ can be analytically found as the
solution of linear ordinary differential equation with constant coefficients.
In practice, we truncate the expansion \eqref{eq_series_waveguide} at $n=50$,
which is sufficient since the convergence is exponential. $\bh$ is easily recovered
by (analytically) differentiating \eqref{eq_series_waveguide}.

First, we consider $\omega = 2\pi$ with a $64 \times 16$ Cartesian mesh and $\CP_1$ elements.
Then, for $\omega = 6\pi$ we use $\CP_3$ elements on a $32 \times 8$ Cartesian mesh.

Figures \ref{figure_waveguide_error} shows the convergence history of the FW and CMCG solver.
CMCG converges significantly faster than FW. In particular, for $\omega = 6\pi$, the
FW solver does not reach convergence within 1000 simulated periods. As in the previous experiment,
CMCG achieves the same accuracy than FS for $\CP_1$ elements, while the error is slightly
increased for $\CP_3$ elements.

\input{figures/waveguide_error}

\subsection{Cavity problem}
\label{section_cavity_experiment}

We next consider an interior problem in a closed cavity $\Omega \eq (0,1)^2$
surrounded by a conducting material. We thus set $\GP \eq \partial \Omega$
and $\GI \eq \emptyset$. We apply a source $j \eq 1$ and set $g \eq 0$.
This problem features resonances at frequencies $\omega_{{\rm r},n,m}^2 \eq (n^2 + m^2) \pi^2$,
for all $n,m \geq 0$, with associated eigenmodes $u_{n,m} \eq \sin(n\pi \bx_1)\sin(m\pi \bx_2)$.
Again, we obtain a semi-analytical solution by truncating the Fourierexpansion.

We examine the behaviour of FW and CMCG when the frequency
$\omega$ is relatively far or close to a resonant frequency $\omega_{\rm r}$.
Hence, for a fixed resonant frequency $\omega_{\rm r}$, we consider a frequency
of the form $\omega_\delta \eq \omega_{\rm r} + \sqrt{2}\pi \delta$
with $\delta = 1/8$ or $1/64$.
We first take $\omega_{\rm r} \eq 3\sqrt{2}\pi$ with $\CP_1$ elements
and a $32 \times 32$ Cartesian mesh. Then, we use $\CP_3$ elements on
an $8 \times 8$ Cartesian mesh for $\omega_{\rm r} \eq 5\sqrt{2}\pi$.

Figures \ref{figure_cavity_error_O00} and \ref{figure_cavity_error_O02}
depict the convergence history of FW and CMCG. The FW algorithm fails to converge
even in the favorable case where $\delta = 1/8$. The CMCG algorithm converges in all cases,
and the convergence rate is only slightly affected for the smaller value of $\delta$.

\input{figures/cavity_error}

\subsection{Dipole source in a trapping medium}

The goal of this experiment is to modelize the electromagnetic field
generate by a dipole source emitting inside a body $G \subset \Omega \eq (-1,1)^2$.
We set $\GP \eq \emptyset$ and $\GI \eq \partial \Omega$.
The permittivity is not constant, and instead, we assume that
\begin{equation*}
\eps(\bx) \eq \left |
\begin{array}{ll}
4 & \text{ if } \bx \in G,
\\
1 & \text{ otherwise},
\end{array}
\right .
\end{equation*}
this choice is made so that $G$ traps rays: Snell's law ensures that rays
crossing the interface with incident angle less that $60^o$ are totally reflected inside the
$G$. We modelize the dipole with $j(\bx) \eq \exp \left (-|\bx-\bc|^2/s^2\right )$
where $s \eq 0.05$ and $\bc \in G$ is the dipole localization. We consider two configurations.
In the first case, the trapping body $G \eq [-0.5,0.5]^2$ is squared, $\bc \eq (0.25,0)$ and
$\omega \eq 10\pi$. In the second case $G \eq \{ \bx \in \mathbb R^2 \; | \; |\bx| < 0.5 \}$
is a disk, $\bc \eq (\sqrt{2}/4,1/2-\sqrt{2}/4)$ and $\omega \eq 20\pi$.
We employ unstructured meshes generated with {\tt GMSH} \cite{geuzaine_remacle_2009a}.
For the square case, we impose a maximum element size $h=0.05$ leading
to a 3636 elements mesh. For the circular trap, the condition $h=0.02$
leads to a 22294 triangles mesh. In both cases, $\CP_3$ elements are used respectively
resulting in 109k and 668k degrees of freedom. Figure \ref{figure_trap_solution} represents
the solutions while Figure \ref{figure_trap_error} shows the behaviour of the error. Again,
CMCG clearly outperforms FW.

\input{figures/trap_solution}

\input{figures/trap_error}

%% file: figures/planewave_error.tex
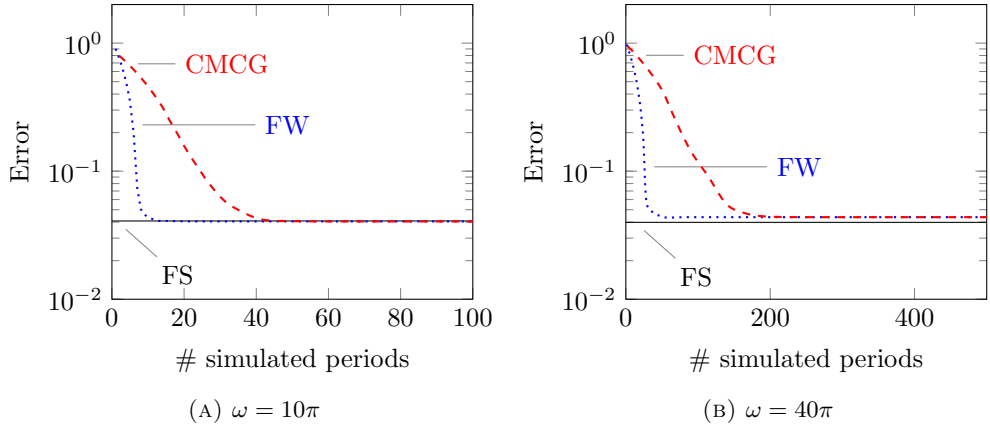
\begin{figure}

\begin{minipage}{.45\linewidth}
\begin{tikzpicture}
\begin{axis}%
[
	width=.95\linewidth,
	xlabel={\# simulated periods},
	ylabel={Error},
	ymode = log,
	xmin=0,
	xmax=100,
	ymax=2,
	ymin=1.e-2
]

\addplot[black,mark=none] coordinates {(0,0.0408405) (100,0.0408405)}
node [pos=.0125,pin={[pin distance=.5cm,pin edge=solid]-45:{FS}}] {};

\addplot[blue,thick,dotted]
table[mark=none,x expr =1+\thisrow{N},y=err]
{figures/data/planewave/F005.00_P045.00_N0032_O00_C00.240/full_wave.txt}
node [pos=.05,pin={[pin distance=1.5cm, pin edge=solid]0:{FW}}] {};

\addplot[red,thick,dashed]
table[mark=none,x expr =2*(1+\thisrow{N}),y=err]
{figures/data/planewave/F005.00_P045.00_N0032_O00_C00.240/cmcg.txt}
node [pos=.0250,pin={[pin distance=.5cm, pin edge=solid]0:{CMCG}}] {};

\end{axis}
\end{tikzpicture}
\subcaption{$\omega = 10\pi$}
\end{minipage}
\begin{minipage}{.45\linewidth}
\begin{tikzpicture}
\begin{axis}%
[
	width=.95\linewidth,
	xlabel={\# simulated periods},
	ylabel={Error},
	ymode = log,
	xmin=0,
	xmax=500,
	ymax=2,
	ymin=1.e-2
]

\addplot[black,mark=none] coordinates {(0,0.0398453) (500,0.0398453)}
node [pos=.025,pin={[pin distance=.5cm, pin edge=solid]-45:{FS}}] {};

\addplot[blue,thick,dotted]
table[mark=none,x expr =1+\thisrow{N},y=err]
{figures/data/planewave/F020.00_P045.00_N0032_O02_C00.120/full_wave.txt}
node [pos=.05,pin={[pin distance=1.5cm, pin edge=solid]0:{FW}}] {};

\addplot[red,thick,dashed]
table[mark=none,x expr =2*(1+\thisrow{N}),y=err]
{figures/data/planewave/F020.00_P045.00_N0032_O02_C00.120/cmcg.txt}
node [pos=.0250,pin={[pin distance=.5cm, pin edge=solid]0:{CMCG}}] {};

\end{axis}
\end{tikzpicture}
\subcaption{$\omega = 40\pi$}
\end{minipage}

\caption{Convergence in the planewave experiment}
\label{figure_planewave_error}
\end{figure}

%% file: figures/waveguide_error.tex
\begin{figure}

\begin{minipage}{.45\linewidth}
\begin{tikzpicture}
\begin{axis}%
[
	width=.90\linewidth,
	xlabel={\# simulated periods},
	ylabel={Error},
	ymode = log,
	xmin=0,
	xmax=1000,
	ymax=2,
	ymin=1.e-3
]

\addplot[black,mark=none] coordinates {(0,0.0139033) (1000,0.0139033)}
node [pos=.025,pin={[pin distance=.5cm,pin edge=solid]-45:{FS}}] {};

\addplot[blue,thick,dotted]
table[mark=none,x expr =1+\thisrow{N},y=err]
{figures/data/waveguide/F001.00_P030.00_N0016_O00_C00.240/full_wave.txt}
node [pos=.65,pin={[pin distance=.5cm,pin edge=solid]45:{FW}}] {};

\addplot[red,thick,dashed]
table[mark=none,x expr =2*(1+\thisrow{N}),y=err]
{figures/data/waveguide/F001.00_P030.00_N0016_O00_C00.240/cmcg.txt}
node [pos=.075,pin={[pin distance=1cm,pin edge=solid]45:{CMCG}}] {};

\end{axis}
\end{tikzpicture}
\subcaption{$\omega = 2\pi$}
\end{minipage}
\begin{minipage}{.45\linewidth}
\begin{tikzpicture}
\begin{axis}%
[
	width=.90\linewidth,
	xlabel={\# simulated periods},
	ylabel={Error},
	ymode = log,
	xmin=0,
	xmax=1000,
	ymax=2,
	ymin=1.e-3
]

\addplot[black,mark=none] coordinates {(0,0.00714692) (1000,0.00714692)}
node [pos=.025,pin={[pin distance=.5cm,pin edge=solid]-45:{FS}}] {};

\addplot[blue,thick,dotted]
table[mark=none,x expr =1+\thisrow{N},y=err]
{figures/data/waveguide/F003.00_P030.00_N0008_O02_C00.120/full_wave.txt}
node [pos=.8,pin={[pin distance=.25cm,pin edge=solid]90:{FW}}] {};

\addplot[red,thick,dashed]
table[mark=none,x expr =2*(1+\thisrow{N}),y=err]
{figures/data/waveguide/F003.00_P030.00_N0008_O02_C00.120/cmcg.txt}
node [pos=.55,pin={[pin distance=.5cm, pin edge=solid]180:{CMCG}}] {};

\end{axis}
\end{tikzpicture}
\subcaption{$\omega = 6\pi$}
\end{minipage}

\caption{Convergence in the half open waveguide experiment}
\label{figure_waveguide_error}
\end{figure}
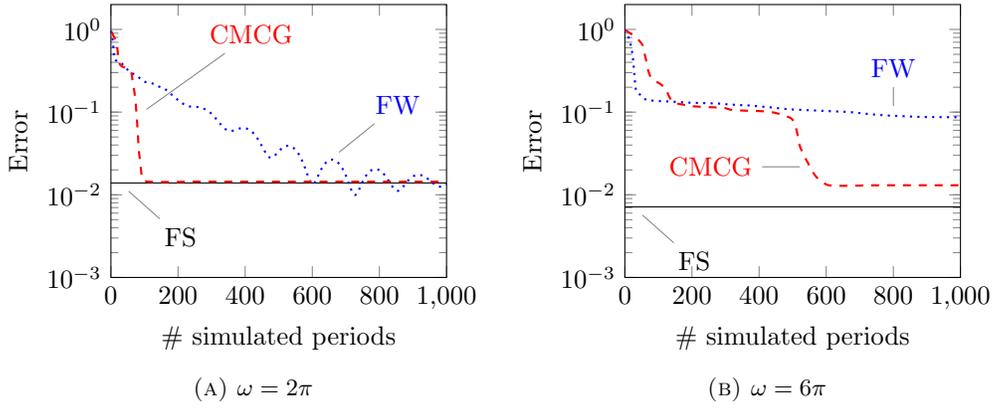

%% file: figures/cavity_error.tex
\begin{figure}
\begin{minipage}{.45\linewidth}
\begin{tikzpicture}
\begin{axis}%
[
	width=.95\linewidth,
	xlabel={\# simulated periods},
	ylabel={Error},
	ymode = log,
	xmin=1,
	xmax=100,
	ymin=1.e-3,
	ymax=2.e-0
]

\addplot[black,thick,solid] coordinates {(1,0.00457162) (100,0.00457162)}
node [pos=.025,pin={[pin distance=.25cm,pin edge=solid]-45:{FS}}] {};

\addplot[blue,thick,dotted]
table[mark=none,x expr =1+\thisrow{N},y=err]
{figures/data/cavity/F02.209707000000_N0032_O00_C00.240/full_wave.txt}
node [pos=.8,pin={[pin distance=.25cm,pin edge=solid]-90:{FW}}] {};

\addplot[red,thick,dashed]
table[mark=none,x expr =2*(1+\thisrow{N}),y=err]
{figures/data/cavity/F02.209707000000_N0032_O00_C00.240/cmcg.txt}
node [pos=.1,pin={[pin distance=.25cm,pin edge=solid]0:{CMCG}}] {};

\end{axis}
\end{tikzpicture}
\subcaption{$\delta = 1/8$}
\end{minipage}
\begin{minipage}{.45\linewidth}
\begin{tikzpicture}
\begin{axis}%
[
	width=.95\linewidth,
	xlabel={\# simulated periods},
	ylabel={Error},
	ymode = log,
	xmin=1,
	xmax=100,
	ymin=1.e-3,
	ymax=2.e-0
]

\addplot[black,thick,solid] coordinates {(1,0.024144) (100,0.024144)}
node [pos=.025,pin={[pin distance=.5cm,pin edge=solid]-45:{FS}}] {};

\addplot[blue,thick,dotted]
table[mark=none,x expr =1+\thisrow{N},y=err]
{figures/data/cavity/F02.132368000000_N0032_O00_C00.240/full_wave.txt}
node [pos=.8,pin={[pin distance=.25cm,pin edge=solid]-90:{FW}}] {};

\addplot[red,thick,dashed]
table[mark=none,x expr =2*(1+\thisrow{N}),y=err]
{figures/data/cavity/F02.132368000000_N0032_O00_C00.240/cmcg.txt}
node [pos=.3,pin={[pin distance=.25cm,pin edge=solid]0:{CMCG}}] {};

\end{axis}
\end{tikzpicture}
\subcaption{$\delta = 1/64$}
\end{minipage}

\caption{Convergence in the cavity experiment: $\omega_{\rm r} = 3\sqrt{2}\pi$}
\label{figure_cavity_error_O00}
\end{figure}

\begin{figure}

\begin{minipage}{.45\linewidth}
\begin{tikzpicture}
\begin{axis}%
[
	width=.95\linewidth,
	xlabel={\# simulated periods},
	ylabel={Error},
	ymode = log,
	xmin=1,
	xmax=100,
	ymin=1.e-3,
	ymax=2.e-0
]

\addplot[black,solid] coordinates {(1,0.00675396) (100,0.00675396)}
node [pos=.025,pin={[pin distance=.5cm,pin edge=solid]-45:{FS}}] {};

\addplot[blue,thick,dotted]
table[mark=none,x expr =1+\thisrow{N},y=err]
{figures/data/cavity/F03.623920000000_N0008_O02_C00.120/full_wave.txt}
node [pos=.8,pin={[pin distance=.25cm,pin edge=solid]-90:{FW}}] {};

\addplot[red,thick,dashed]
table[mark=none,x expr =2*(1+\thisrow{N}),y=err]
{figures/data/cavity/F03.623920000000_N0008_O02_C00.120/cmcg.txt}
node [pos=.2,pin={[pin distance=.25cm,pin edge=solid]0:{CMCG}}] {};

\end{axis}
\end{tikzpicture}
\subcaption{$\delta = 1/8$}
\end{minipage}
\begin{minipage}{.45\linewidth}
\begin{tikzpicture}
\begin{axis}%
[
	width=.95\linewidth,
	xlabel={\# simulated periods},
	ylabel={Error},
	ymode = log,
	xmin=1,
	xmax=100,
	ymin=1.e-3,
	ymax=2.e-0
]

\addplot[black,solid] coordinates {(1,0.0284405) (100,0.0284405)}
node [pos=.025,pin={[pin distance=.5cm,pin edge=solid]-45:{FS}}] {};

\addplot[blue,thick,dotted]
table[mark=none,x expr =1+\thisrow{N},y=err]
{figures/data/cavity/F03.546581000000_N0008_O02_C00.120/full_wave.txt}
node [pos=.8,pin={[pin distance=.25cm,pin edge=solid]-90:{FW}}] {};

\addplot[red,thick,dashed]
table[mark=none,x expr =2*(1+\thisrow{N}),y=err]
{figures/data/cavity/F03.546581000000_N0008_O02_C00.120/cmcg.txt}
node [pos=.3,pin={[pin distance=.25cm,pin edge=solid]0:{CMCG}}] {};

\end{axis}
\end{tikzpicture}
\subcaption{$\delta = 1/64$}
\end{minipage}

\caption{Convergence in the cavity experiment: $\omega_{\rm r} = 5\sqrt{2}\pi$}
\label{figure_cavity_error_O02}
\end{figure}

%% file: figures/trap_solution.tex
\begin{figure}

\begin{minipage}{.50\linewidth}
\begin{tikzpicture}

\node[inner sep=0pt,anchor=south west] at (0,0) {\includegraphics[width=5cm,height=5cm]%
{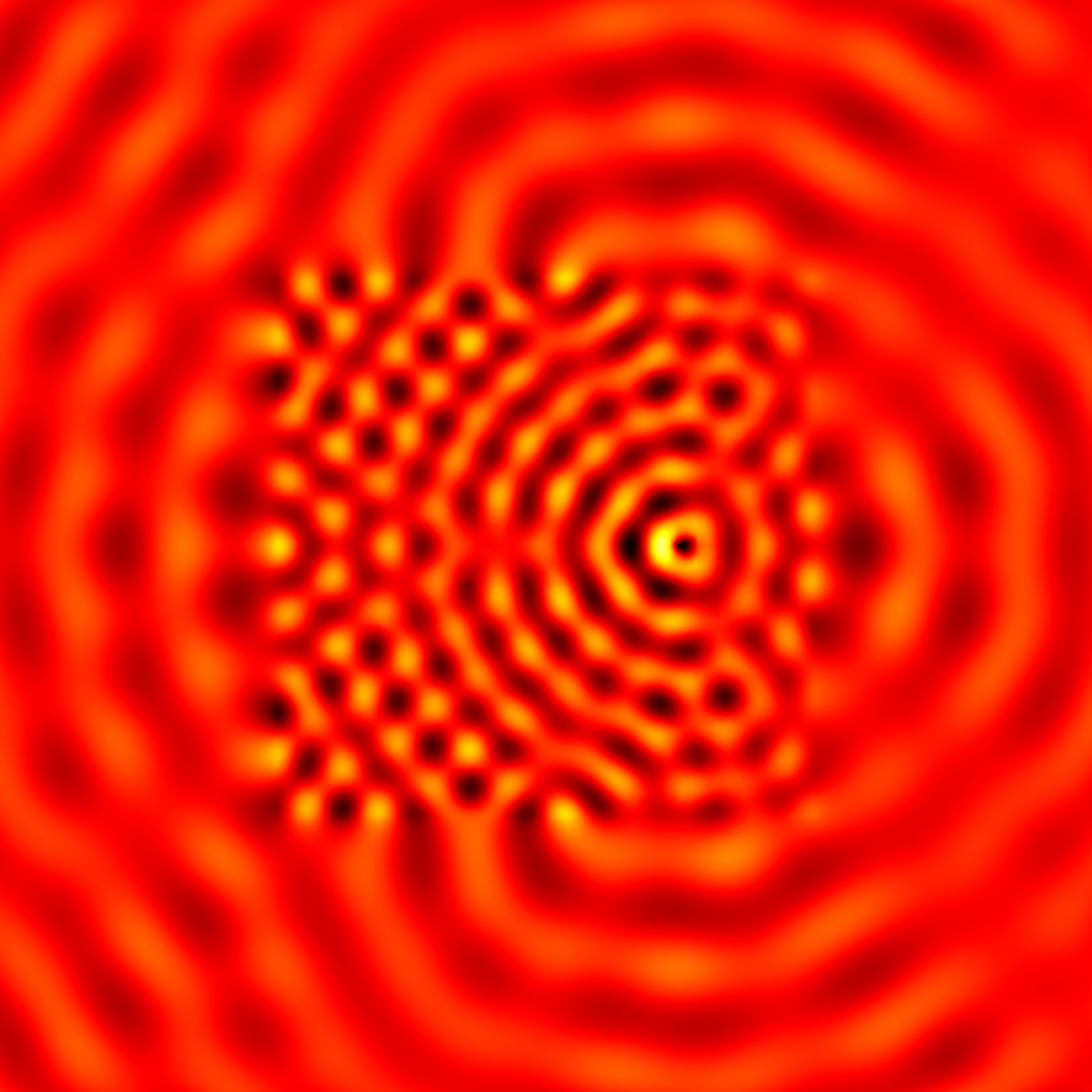}};

\node[anchor=north] at (0.0,0) {-1};
\node[anchor=north] at (5.0,0) { 1};

\node[anchor=east]  at (0,0.0) {-1};
\node[anchor=east]  at (0,5.0) { 1};

\draw[thick] (0,0  ) rectangle (5,5);
\end{tikzpicture}
\end{minipage}
\begin{minipage}{.40\linewidth}
\begin{tikzpicture}

\node[inner sep=0pt,anchor=south west] at (0,0) {\includegraphics[width=5cm,height=5cm]%
{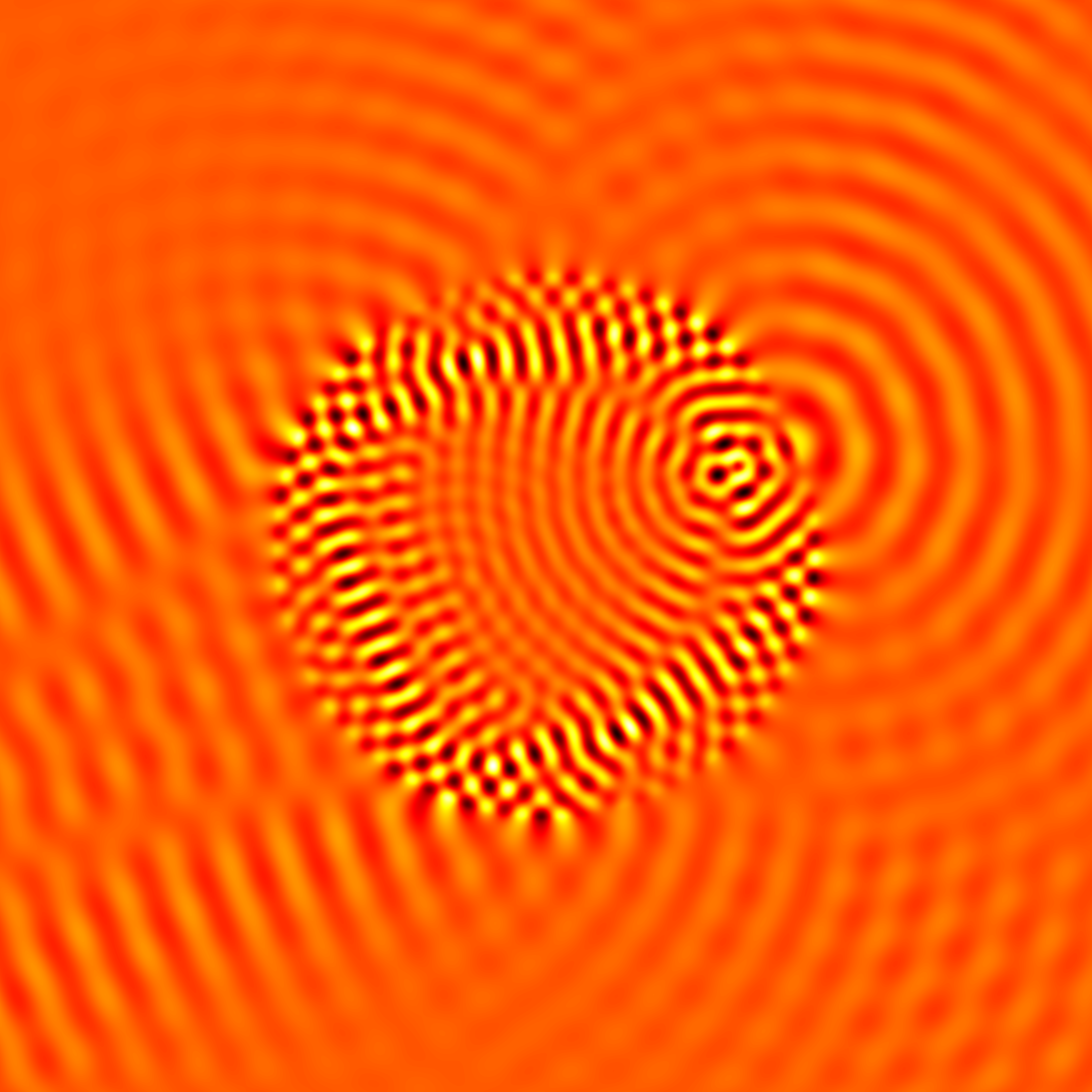}};

\node[anchor=north] at (0.0,0) {-1};
\node[anchor=north] at (5.0,0) { 1};

\node[anchor=east]  at (0,0.0) {-1};
\node[anchor=east]  at (0,5.0) { 1};

\draw[thick] (0,0  ) rectangle (5,5);

\node[inner sep=0pt,anchor=south west] at (5.5,0) {\includegraphics[width=.5cm,height=5cm]%
{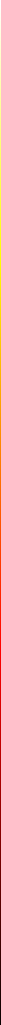}};

\draw[thick] (5.5,0) rectangle (6,5);

\node[anchor=west] at (6,0.0) {$-1.5 \times 10^{-3}$};
\node[anchor=west] at (6,5.0) {$ 2.5 \times 10^{-3}$};
\end{tikzpicture}
\end{minipage}

\caption{Imaginary part of the electric field in the square (left) and circular (right) traps}
\label{figure_trap_solution}
\end{figure}
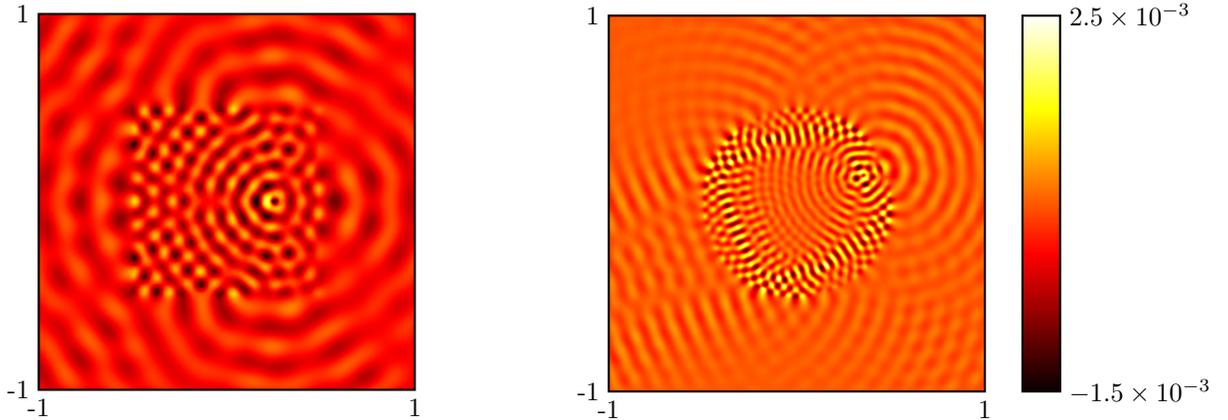

%% file: figures/trap_error.tex
\begin{figure}

\begin{minipage}{.45\linewidth}
\begin{tikzpicture}
\begin{axis}%
[
	width=.90\linewidth,
	xlabel={\# simulated periods},
	ylabel={Error},
	ymode = log,
	xmin=0,
	xmax=1000,
	ymax=2,
	ymin=5.e-3
]

\addplot[blue,thick,dotted]
table[mark=none,x expr =1+\thisrow{N},y=err]
{figures/data/square_trap/F005.00_P000.00_H0.05_O02_C00.120/full_wave.txt}
node [pos=.5,pin={[pin edge=solid]45:{FW}}] {};

\addplot[red,thick,dashed]
table[mark=none,x expr =2*(1+\thisrow{N}),y=err]
{figures/data/square_trap/F005.00_P000.00_H0.05_O02_C00.120/cmcg.txt}
node [pos=.1,pin={[pin edge=solid]0:{CMCG}}] {};

\end{axis}
\end{tikzpicture}
\end{minipage}
\begin{minipage}{.45\linewidth}
\begin{tikzpicture}
\begin{axis}%
[
	width=.90\linewidth,
	xlabel={\# simulated periods},
	ylabel={Error},
	ymode = log,
	xmin=0,
	xmax=1000,
	ymax=2,
	ymin=5.e-3
]

\addplot[blue,thick,dotted]
table[mark=none,x expr =1+\thisrow{N},y=err]
{figures/data/circle_trap/F010.00_P000.00_H0.02_O02_C00.120/full_wave.txt}
node [pos=.8,pin={[pin edge=solid]-90:{FW}}] {};

\addplot[red,thick,dashed]
table[mark=none,x expr =2*(1+\thisrow{N}),y=err]
{figures/data/circle_trap/F010.00_P000.00_H0.02_O02_C00.120/cmcg.txt}
node [pos=.5,pin={[pin edge=solid]180:{CMCG}}] {};

\end{axis}
\end{tikzpicture}
\end{minipage}

\caption{Convergence in the square (left) and circular (right) trap experiments}
\label{figure_trap_error}
\end{figure}
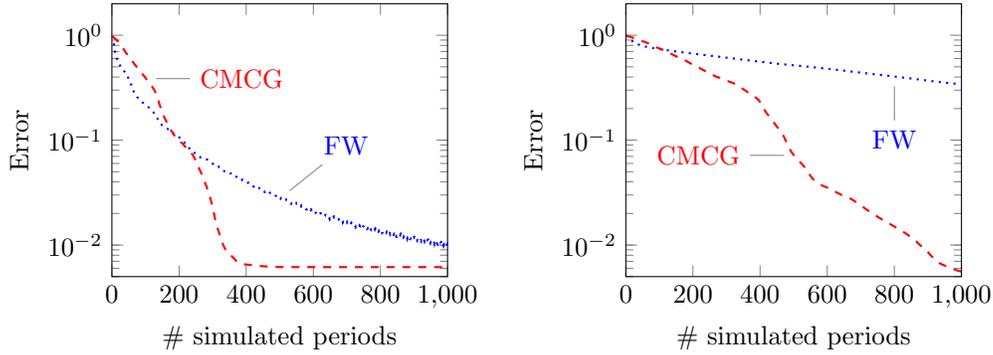

%% file: 7-conclusion.tex
\section{Conclusion}
\label{section_conclusion}

We propose a controllability method (CM) to solve Maxwell's equations in the frequency-domain
in their first-order formulation. By minimizing a quadratic
cost functional $J$ using a conjugate gradient iteration (CG), the CMCG method determines
a  time-periodic solution in the time-domain. At each CG iteration, the gradient $J'$ is computed
simply by running a time-domain solver forward and backward over one period, without
the need for solving any additional linear system. Hence, our CMCG algorithm automatically inherits
the parallelism, scalability, and low memory footprint of the underlying DG
time-domain solver. The full CMCG Algorithm 2.3  is listed in Section \ref{sec_CMCG_algo}.

In general, there exist several time-periodic solutions to Maxwell's equations,
distinct from the desired  time-harmonic solution,
so that the minimizer of $J$ may not be unique. To remove those spurious modes and thus extract
the time-harmonic solution from any minimizer, we apply a cheap filtering operator computed
``on the fly'' as a final post-processing step. In Theorem \ref{theorem_filtering},
we establish that $J$ combined with the filtering operator is
strongly convex in an appropriate sense, which ensures the convergence of the CMCG method
to the desired time-harmonic solution from any initial guess. In Section \ref{section_filtering},
we also show that nearly periodic solutions already provide good approximations to
the time-harmonic solution after filtering. Hence, by monitoring the misfit, the CG iteration
may be stopped as soon as the desired accuracy has been reached.

Comparison with a direct frequency-domain solver shows that the additional
error due to time discretization is hardly visible for the low-order
$\CP_1$-RK2 discretization and very small for the higher order $\CP_3$-RK4 discretization.
In these numerical experiments, we also compare the CMCG method against the limiting amplitude
principle, where one simply lets the time-domain solver run until the time-harmonic regime
is reached. For simple plane wave propagation, the
limiting amplitude principle in fact slightly outperforms CMCG. For all other examples however,
CMCG significantly outperforms the limiting amplitude approach.
For the cavity experiment in Section
\ref{section_cavity_experiment}, in particular, the convergence of
CMCG is hardly affected by the trapping geometry,
whereas the limiting amplitude principle utterly fails.

Our CMCG method is non-intrusive and easily integrated
into any existing time-domain code. It is
not limited to DG discretizations; thus, we expect similar performance
using solvers based on finite differences \cite{taflove_hagness_2005a,yee_1966a}.
Although we have only used simple first-order Silver-Müller absorbing boundary conditions
in our computations, the CMCG approach immediately extends to other more accurate absorbing
conditions or perfectly matched layers \cite{tang_2020a}.
In the presence of complex geometry and local mesh refinement,  local time-stepping methods
permit to overcome the stringent local CFL stability condition without sacrificing explicitness
\cite{grote_mehlin_mitkova_2015a,grote_tang_2019a}. The CMCG approach can also compute solutions
for multiple frequencies in ``one shot'', that is at the cost of a single solve, as proposed in
\cite{tang_2020a}.